\documentclass[preprint,12pt,compress]{elsarticle}

\usepackage{amssymb}
\usepackage{amsmath}
 \usepackage{bbm}
\usepackage[a4paper, body={14.8cm,23cm}]{geometry}

\usepackage{amsthm}
\usepackage{enumitem}
\usepackage{fancyhdr}
\usepackage{hyperref}
\hypersetup{%
  pdftitle={Existence, uniqueness and approximation for $L^p$ solutions of reflected BSDEs under weaker assumptions},%
  pdfsubject={L^p$ solutions of reflected BSDEs},%
  pdfauthor={ShengJun FAN},%
  pdfkeywords={BSDEs, Stability theorem, comparison theorem, Existence and uniqueness, $L^p$ solution},%
  pdfstartview=FitH,
  CJKbookmarks=true,%
  bookmarksnumbered=true,%
  bookmarksopen=true,%
  colorlinks=true, linkcolor=blue, urlcolor=blue, citecolor=blue, %
}

\begin{document}
\journal{arXiv}
\parindent 15pt
\parskip 5pt
\renewcommand{\arraystretch}{2}

 \newcommand{\eps}{\varepsilon}
 \newcommand{\lam}{\lambda}
 \newcommand{\To}{\rightarrow}
 \newcommand{\as}{{\rm d}\mathbb{P}\times{\rm d}t-a.e.}
 \newcommand{\ass}{{\rm d}\mathbb{P}\times{\rm d}s-a.e.}
 \newcommand{\ps}{\mathbb{P}-a.s.}
 \newcommand{\jf}{\int_t^T}
 \newcommand{\tim}{\times}

 \newcommand{\F}{\mathcal{F}}
 \newcommand{\E}{\mathbb{E}}
 \newcommand{\N}{\mathbb{N}}

 \newcommand{\s}{\mathcal{S}}
 \newcommand{\M}{{\rm M}}
 \newcommand{\hcal}{\mathcal{H}}
 \newcommand{\vcal}{\mathcal{V}}
 \newcommand{\mcal}{\mathcal{M}}

 \newcommand{\T}{[0,T]}
 \newcommand{\Lp}{{\mathbb L}^p(\F_T)}

 \newcommand{\R}{{\mathbb R}}
 \newcommand{\Q}{{\mathbb Q}}
 \newcommand{\RE}{\forall}

\newcommand {\Lim}{\lim\limits_{n\rightarrow \infty}}
\newcommand {\Dis}{\displaystyle}

\begin{frontmatter}
\title {Existence, uniqueness and approximation for $L^p$ solutions of reflected BSDEs under weaker assumptions\tnoteref{fund}}
\tnotetext[fund]{Supported by the National Natural Science Foundation of China (No. 11371362), the China Postdoctoral Science Foundation (No. 2014T70386 and 2013M530173), the Qing Lan Project and the Fundamental Research Funds for the Central Universities (No. 2013RC20)\vspace{0.2cm}.}

\author{ShengJun FAN}
\ead{f$\_$s$\_$j@126.com}

\address{College of Sciences, China University of Mining and Technology, Xuzhou 221116, PR China\vspace{-0.8cm}}

\begin{abstract}
We put forward and prove several existence and uniqueness results for $L^p\ (p>1)$ solutions of reflected BSDEs with continuous barriers and generators satisfying a one-sided Osgood condition together with a general growth condition in $y$ and a uniform continuity condition or a linear growth condition in $z$. A necessary and sufficient condition with respect to the growth of barrier is also explored to ensure the existence of a solution. And, we show that the solutions may be approximated by the penalization method and by some sequences of solutions of reflected BSDEs. Our results improve considerably some known works.\vspace{0.1cm}
\end{abstract}

\begin{keyword}
Reflected backward stochastic differential equation\sep $L^p$ solutions\sep Comparison theorem \sep One-sided Osgood condition \sep  Uniform continuity condition\vspace{0.3cm}

\MSC[2010] 60H10
\end{keyword}
\end{frontmatter}\vspace{-0.4cm}


\section{Introduction}

Let $(\Omega,\F,\mathbb{P})$ be a completed probability space carrying a standard $d$-dimensional Brownian motion $(B_t)_{t\geq 0}$, and $(\F_t)_{t\geq 0}$ the completed $\sigma$-algebra filtration generated by $(B_t)_{t\geq 0}$. Assume that $T>0$ is a real number and $\F=\F_T$. In this paper we are given an $(\F_t)_{t\geq 0}$-progressively measurable continuous process $(L_t)_{t\in\T}$, an $\F_T$-measurable random variable $\xi$ such that $\xi\geq L_T$, a random function
$$g(\omega,t,y,z):\Omega\tim \T\tim {\R}\tim
{\R}^{d}\longmapsto {\R}$$
such that $g(\cdot,y,z)$ is $(\F_t)$-progressively measurable for each $(y,z)$, and an $(\F_t)_{t\geq 0}$-progressively measurable continuous process $(V_t)_{t\in\T}$ with finite variation. By a solution of the reflected backward stochastic differential equation (reflected BSDE or directly RBSDE for short) with terminal time $T$, terminal value $\xi$, generator $g+dV$ and barrier $L_\cdot$ we understand a triple $(Y_t,Z_t,K_t)_{t\in\T}$ of $(\F_t)_{t\geq 0}$-progressively measurable processes such that $t\mapsto |g(t,Y_t,Z_t)|$ belongs to ${\rm L}^1(0,T)$, $t\mapsto |Z_t|$ belongs to ${\rm L}^2(0,T)$ and\vspace{0.1cm}
\begin{equation}
\hspace*{-0.5cm}\left\{
\begin{array}{l}
\Dis Y_t=\xi+\int_t^Tg(s,Y_s,Z_s){\rm d}s+\int_t^T{\rm d}V_s+\int_t^T{\rm d}K_s-\int_t^TZ_s {\rm d}B_s,\ \   t\in\T,\\
\Dis Y_t\geq L_t,\ \ t\in\T,\\
\Dis K\ {\rm is\ nondecreasing,\ continuous},\  K_0=0,\ \ \int_0^T (Y_t-L_t){\rm d}K_t
=0.
\end{array}
\right.
\end{equation}
This equation is usually denoted by RBSDE $(\xi,g+dV,L)$. The second condition in (1) says that the first component $Y_\cdot$ of the solution is forced to stay above $L_\cdot$. The role of $K_\cdot$ is to push $Y_\cdot$ upwards in order to keep it above $L_\cdot$ in a minimal way, which means that the third condition in (1) is satisfied. Note that the usual BSDEs may be considered as a special case of RBSDEs with $L_\cdot\equiv -\infty$ (and  $K_\cdot\equiv 0$).

Nonlinear BSDEs were initially put forward in 1990 by \citet{Par90}, which proved an existence and uniqueness result for square-integrable solutions of BSDEs with generators satisfying the Lipschitz condition in $(y,z)$ where the data $\xi,g(\cdot,0,0)$ are square-integrable and $V_\cdot\equiv 0$. In \citet{El97b} and \citet{Peng99}, the authors further investigated BSDEs with $V_\cdot$ being not zero. As a generation of the notion of nonlinear BSDEs, \citet{El97} first introduced nonlinear RBSDEs in 1997 and proved the existence and uniqueness for square-integrable solutions of RBSDEs with generators satisfying the Lipschitz condition in $(y,z)$ where the data $\xi,g(\cdot,0,0)$ and $\sup_{t\in\T}|L_t|$ are all square-integrable and $V_\cdot\equiv 0$. Recently, \citet{Kli13} further considered RBSDEs with $V_\cdot$ being not zero and with discontinuous barriers. BSDEs and RBSDEs have attracted more and more interests, and due to the closely connections with many problems, they have gradually become a very useful and efficient tool in different mathematical fields including mathematical finance, game theory, optimal switching problem, partial differential equations and others (see, e.g., \cite{Bay12,Bay14,El97,El97a,El97b,
Ham02,Ham00,Ham10,Hu08,Jia10,Par99,Peng99,Peng05,Peng10}).

The assumptions on the data in \cite{El97,Par90} are sometimes too strong for many interesting applications. Therefore many attempts have been made to prove the existence and uniqueness of solutions of BSDEs or RBSDEs under less restrictive assumptions. For example, many papers were devoted to relaxing the continuity assumption for the barrier $L_\cdot$ of RBSDEs, see \cite{Bay14,Ham02,Kli13,Lep05a,Peng05}; many papers aimed to solving BSDEs or RBSDEs with data that are not square-integrable but only in $L^p\ (p>1)$ or $L^1$, see \cite{Ama09,Bay14,Bri03,Bri07,El97b,
Fan15,Ham12,Izu13,Kli12,Kli13,Ma13,Roz12}; and more papers were interested in weakening the linear growth and Lipschitz-continuity of the generator $g$ with respect to $(y,z)$, see \cite{Bri03,Bri07,El97b,Fan15,Fan10,Fan12,Fan13,
Hu14,Izu13,Jia08,Jia10,Kob02,Ma13,Mao95,Par99} for BSDEs and \cite{Ama09,Bay12,Bay14,Ham10,Hua13,Jia08a,
Kli12,Kli13,Lep05,Mat97,Roz12,Xu08} for RBSDEs.

In the present paper we focus attention on solving RBSDE (1) with $L^p\ (p>1)$ data and continuous barriers under some weaker assumptions of the generator $g$ with respect to $(y,z)$. RBSDEs with $L^1$ data and discontinuous barriers under weaker assumptions will be our object of research in the near future. Here, we would like to mention some known results related closely to our work. Firstly, in \citet{Bri03}, the authors established the existence and uniqueness of $L^p$ solutions of BSDEs with $L^p\ (p\geq 1)$ data under assumptions that the generator $g$ satisfies a monotonicity condition together with a general growth condition in $y$ and the Lipschitz condition in $z$ (see respectively (H1s), (H3) and (H2s) in Section 2 for details). Recently, this result was extended to RBSDE (1) with $V_\cdot\equiv 0$ by \citet{Lep05} and \citet{Roz12}, where some additional assumptions relating the growth of generator $g$ with that of barrier $L_\cdot$ were put forward. Under the same assumptions with respect to the generator $g$ as in \cite{Bri03,Roz12}, \citet{Kli12} further explored a necessary and sufficient condition with respect to the growth of barrier $L_\cdot$ (see (H6) in Section 2 for details) to ensure the existence and uniqueness of $L^p$ solutions of RBSDE (1) with $L^p\ (p\geq 1)$ data, continuous barrier $L_\cdot$ and $V_\cdot\equiv 0$. And, by establishing and applying a general monotonic limit theorem of BSDEs, \citet{Kli13} investigated RBSDEs with two irregular reflecting barriers.

Furthermore, in \citet{Mat97} and \citet{Xu08} the authors obtained respectively an existence result of square-integrable solutions for RBSDE (1) with $L^2$ data and $V_\cdot\equiv 0$ where the generator $g$ satisfies a linear growth condition in $z$ instead of the Lipschitz condition. \citet{Jia08a} further proved a uniqueness result where the generator $g$ satisfies a uniform continuity condition in $z$ ( see (H2) in Section 2 for details). With respect to this condition, we also refer to \cite{Fan10,Fan12,Jia08,Jia10}.

On the other hand, \citet{Fan12} put forward a kind of one-sided Osgood condition in $y$ of the generator $g$ (see (H1) in Section 2 for details), which is weaker than not only Mao's condition used in \citet{Mao95} and the Osgood condition put forward in \citet{Fan13} but also the monotonicity condition (see (H1s) in Section 2) applied in \cite{Bay14,Bri03,Kli12,Kli13,
Lep05,Roz12,Xu08}. Under (H1) and (H2), a comparison theorem of square-integral solutions of BSDEs was proved in \cite{Fan12}, whose a direct corollary is the uniqueness of solutions. Furthermore, the existence and uniqueness result for $L^p$ solutions of BSDEs with $L^p\ (p>1)$ data obtained in \cite{Bri03} has been extended by \citet{Fan15}, where the generator $g$ satisfies (H1), (H3) and (H2s). Then, the following questions are naturally asked:
\begin{itemize}
\item Can we establish a comparison theorem for $L^p$ solutions of RBSDEs with $L^p\ (p>1)$ data under assumptions (H1) and (H2)?

\item Can we establish an existence result for $L^p$ solutions of RBSDE (1) with $L^p\ (p>1)$ data under some appropriate conditions with respect to the growth of barrier $L_\cdot$ if the generator $g$ only satisfies (H1), (H3) and (H2) or (H2w), a linear growth condition in $z$ (see Section 2 for details)?

\item Can we give a necessary and sufficient condition with respect to the growth of barrier $L_\cdot$ to ensure the existence of $L^p$ solutions of RBSDE (1) with $L^p\ (p>1)$ data under (H1), (H3) and (H2) or (H2w)?

\item Does the sequence of $L^p$ solutions of usual penalization equation for RBSDE (1) with $L^p\ (p>1)$ data still converge under (H1), (H3) and (H2) or (H2w)?
\end{itemize}
The present paper gives positive answers for all these questions. It should be mentioned that our results improve considerably some works mentioned before and that many technical results in our work, including some a priori estimates of BSDEs, the convergence of sequence of $L^p$ solutions for penalization and approximation equations of RBSDE (1) with $L^p\ (p>1)$ data and the comparison theorem for $L^p$ solutions of RBSDEs, are all respectively established under some very general and elementary conditions, for example, the assumptions (H1), (H2), (HH) and (A) in Section 2 and the assumptions (B) and (C) in Section 3.

The remainder of this paper is organized as follows. Section 2 contains some notation and hypotheses which will be used later. Section 3 is devoted to establishing several a priori estimates on $L^p$ solutions of BSDEs with $L^p\ (p>0)$ data as well as on some sequences of $L^p$ solutions of BSDEs with $L^p\ (p>1)$ data, which will play important roles in the proof of our main results. The convergence of sequence of $L^p$ solutions for penalization and approximation equations of RBSDE (1) with $L^p\ (p>1)$ data under assumption (HH) and some very elementary conditions, and a comparison theorem for $L^p$ solutions of RBSDEs with $L^p\ (p>1)$ data under assumptions (H1) and (H2) are put forward and proved in Section 4. Based on these results, Section 5 focus on establishing some existence, uniqueness and approximation results on $L^p$ solutions of BSDEs with $L^p\ (p>1)$ data (Section 5.1) and RBSDEs with $L^p\ (p>1)$ data (Section 5.2) under weaker assumptions, which answers those questions put forward before.

We note that the work of this paper (even for Theorems 1-2 on non-reflected BSDEs) improves considerably some corresponding known results including those obtained in \citet{Bri07,El97b,Fan15,Fan12,Ham12,Kli12,Lep05,Roz12} and \citet{Xu08} (see, for example, Remark 5 in Section 4 and Remarks 6 and 8 in Section 5.2 for more details).

\section{Notation and hypotheses}

In the whole paper we fix a real number $T>0$ and a positive integer $d$, and let $\R_+:=[0,+\infty)$, $a^+:=\max\{a,0\}$ and $a^-:=(-a)^+$ for any real number $a$. Let $\mathbbm{1}_{A}$ represent the indicator function of a set $A$, and ${\rm sgn}(x)$ the sign of a real number $x$. The Euclidean norm of a vector $z\in \R^{d}$ is denoted by $|z|$.

Let $(\Omega,\F,\mathbb{P})$ be a completed probability space carrying a standard $d$-dimensional Brownian motion $(B_t)_{t\geq 0}$, and $(\F_t)_{t\geq 0}$ the completed $\sigma$-algebra filtration generated by $(B_t)_{t\geq 0}$ and assume that $\F=\F_T$. In the whole paper all equalities and inequalities between random elements are understood to hold $\ps$

For $p>0$, denote by $\Lp$ the set of all $\F_T$-measurable random variables $\xi$ such that
$$\|\xi\|_{\mathbb{L}^p}:=\left(\E[|\xi|^p]\right)^{1\wedge 1/p}<+\infty,$$
and define the following spaces of processes or functions:
\begin{itemize}
\item [$\s$---]  the set of all continuous $(\F_t)$-progressively measurable processes;
\item [$\s^p$---]  the set of all processes $Y_\cdot\in \s$ such that
     $$\|Y\|_{{\s}^p}:=\left(\E[\sup_{t\in\T} |Y_t|^p]\right)^{1\wedge 1/p}<+\infty;$$
\item [$\M$---]  the set of all $(\F_t)$-progressively measurable processes $Z_\cdot$ such that
    $$\mathbb{P}\left(\int_0^T|Z_t|^2{\rm d}t<+\infty\right)=1;$$
\item [$\M^p$---]  the set of all processes $Z_\cdot\in \M$ such that
    $$
    \|Z\|_{\M^p}:=\left\{ \E\left[\left(\int_0^T |Z_t|^2{\rm d}t\right)^{p/2}\right] \right\}^{1\wedge 1/p}<+\infty;
    $$
\item [$\hcal$---]  the set of all $(\F_t)$-progressively measurable processes $X_\cdot$ such that
    $$\mathbb{P}\left(\int_0^T|X_t|{\rm d}t<+\infty\right)=1;$$
\item [$\hcal^p$---]  the set of all processes $X_\cdot\in \hcal$ such that
    $$
    \|X\|_{\hcal^p}:=\left\{ \E\left[\left(\int_0^T |X_t|{\rm d}t\right)^p\right] \right\}^{1\wedge 1/p}<+\infty;
    $$
\item [$\mcal$---]  the set of all continuous local $(\F_t)$-martingales.
\item [$\mcal^p$---]  the set of all martingales $M_\cdot\in \mcal$ such that $\E\left[\left(\langle M\rangle_T\right)^{p/2}\right]<+\infty$;
\item [$\vcal$---]  the set of all continuous $(\F_t)$-progressively measurable processes of finite variation;
\item [$\vcal^p$---] the set of all processes $V_\cdot\in \vcal$ such that $\E\left[|V|^p_T\right]<+\infty$;
\item [$\vcal^+$---]  the set of all continuous $(\F_t)$-progressively measurable increasing processes;
\item [$\vcal^{+,p}$---] the set of all processes $V_\cdot\in \vcal^+$ such that $\E\left[|V|^p_T\right]<+\infty$;
\item [${\bf S}$---]  the set of nonnegative functions $\psi_t(\omega,r):\Omega\times \T\times \R_+\To \R_+$
satisfying the following two conditions:
\begin{itemize}
\item $\as$, the function $r\mapsto \psi_t(\omega,r)$ is increasing and $\psi_t(\omega,0)=0$;
\item for each $r\geq 0$, $\psi_\cdot(r)\in \hcal$.\vspace{0.1cm}
\end{itemize}
\end{itemize}

For $p>0$, we introduce the following hypotheses.
\begin{itemize}
\item [{\bf (H1)}] $g$ satisfies the one-sided Osgood condition in $y$, i.e., there exists a nondecreasing and concave function $\rho(\cdot):\R_+\mapsto \R_+$ with $\rho(0)=0$, $\rho(u)>0$ for $u>0$ and $\int_{0^+} {{\rm d}u\over \rho(u)}=+\infty$ such that $\as$, $\RE\ y_1,y_2\in \R,z\in\R^{d}$,
$$
(g(\omega,t,y_1,z)-g(\omega,t,y_2,z)){\rm sgn}(y_1-y_2)\leq \rho(|y_1-y_2|).
$$

\item [{\bf (H1s)}] $g$ satisfies the monotonicity condition in $y$, i.e., there exists a constant $\mu\in \R$ such that $\as$, $\RE\ y_1,y_2\in \R,z\in\R^{d}$,
$$
(g(\omega,t,y_1,z)-g(\omega,t,y_2,z)){\rm sgn}(y_1-y_2)\leq \mu |y_1-y_2|.
$$

\item [{\bf (H2)}] $g$ satisfies the uniform continuity condition in $z$, i.e., there exists a nondecreasing and continuous function $\phi(\cdot):\R_+\mapsto \R_+$ with $\phi(0)=0$ such that $\as$, $\RE\ y\in\R, z_1,z_2\in\R^{d}$,
$$
|g(\omega,t,y,z_1)-g(\omega,t,y,z_2)|\leq \phi(|z_1-z_2|).
$$

\item [{\bf (H2s)}] $g$ satisfies the Lipschitz condition in $z$, i.e., there exists a nonnegative constant $\lambda$ such that $\as$, $\RE\ y\in\R, z_1,z_2\in\R^{d}$,
$$
|g(\omega,t,y,z_1)-g(\omega,t,y,z_2)|\leq \lambda |z_1-z_2|.
$$

\item [{\bf (H2w)}] $g$ has a stronger linear growth in $z$, i.e., there exists two constants $\mu, \lambda\geq 0$ and a nonnegative process $f_\cdot\in \hcal^p$ such that $\as$, $\RE\ y\in\R, z\in\R^{d}$,
$$
|g(\omega,t,y,z)-g(\omega,t,y,0)|\leq f_t(\omega)+\mu |y|+\lambda |z|.
$$

\item [{\bf (H3)}]\ $g$ has a general growth in $y$, i.e, $\RE r>0, \ \varphi_\cdot(r):=\sup\limits_{|y|\leq r}|g(\cdot,y,0)-g(\cdot,0,0)|$ belongs to the space $\hcal$. And, $g(\cdot,0,0)\in\hcal^p$.

\item [{\bf (H3s)}] $g$ has a linear growth in $y$, i.e., there exists a constant $\mu\geq 0$ and a nonnegative process $f_\cdot\in \hcal^p$ such that $\as$, $\RE\ y\in \R$, $|g(\omega,t,y,0)|\leq f_t(\omega)+ \mu |y|$.

\item [{\bf (H4)}] $g$ is continuous in $(y,z)$, i.e., $\as$, $g(\omega,t,\cdot,\cdot)$ is continuous.

\item [{\bf (H4s)}] $g$ is stronger continuous in $(y,z)$, i.e., $\as$, $\RE\ y\in \R,\ g(\omega,t,y,\cdot)$ is continuous, and $g(\omega,t,\cdot,z)$ is continuous uniformly with respect to $z$.

\item [{\bf (H4w)}] $\as$, $\RE\ z\in {\R^{d}},\ \ g(\omega,t,\cdot,z)$ is continuous.

\item [{\bf (H5)}] $\xi\in \Lp$, $L_\cdot\in \s$ and $L_T\leq \xi$.

\item [{\bf (H6)}] There exists a semi-martingale $X_\cdot\in \M^p+\vcal^p$ such that $g(\cdot,X_\cdot,0)\in \hcal^p$ and for each $t\in \T$, $X_t\geq L_t$.

\item [{\bf (HH)}] $g$ has a certain general growth in $(y,z)$, i.e., there exists a constant $\lambda\geq 0$, a nonnegative process $f_\cdot\in\hcal^p$ and a nonnegative function $\psi_\cdot(r)\in {\bf S}$ such that $\as$, $\RE\ y\in\R,\ z\in\R^{d}$,
$$
|g(\omega,t,y,z)|\leq f_t(\omega)+ \psi_t(\omega,|y|)+\lambda |z|.
$$

\item [{\bf (A)}] There exists two nonnegative constants $\bar\mu$ and $\bar\lambda$ such that $\as$,
$$g(\omega,t,y,z){\rm sgn}(y)\leq \bar f_t(\omega)+\bar\mu |y|+\bar\lambda|z|,\ \ \RE\ y\in\R,\ z\in\R^{d},$$
where $\bar f_t$ is a nonnegative process belonging to $\hcal$.
\end{itemize}

{\bf Remark 1}\ Without loss of generality, we will always assume that the functions $\rho(\cdot)$ and $\phi(\cdot)$ defined respectively in (H1) and (H2) are of linear growth, i.e., there exists a constant $A>0$ such that
$$\RE\ x\in \R_+,\ \ \rho(x)\leq A(x+1)\ \ {\rm and}\ \ \phi(x)\leq A(x+1).$$

{\bf Remark 2}\ It is not very hard to verify the following statements hold:
\begin{itemize}
\item [(i)] (H1s)$\Rightarrow$(H1);\ \ (H2s)$\Rightarrow$(H2)$\Rightarrow$(H2w);
    \ \ (H3s)$\Rightarrow$(H3);\ \ (H4s)$\Rightarrow$(H4)$\Rightarrow$(H4w);
\item [(ii)] (H2w)+(H3)$\Rightarrow$(HH)$\Rightarrow$(H3); If (H2) holds, then (H4w)$\Leftrightarrow$(H4);
\item [(iii)] (H6)$\Rightarrow$ $L_\cdot^+\in \s^p$; If $L_\cdot^+\in\s^p$ and
    $$\left(g(t,\sup\limits_{s\in [0,t]}L_s^+,0)\right)_{t\in\T}\in\hcal^p,$$
    then (H6) holds; If (H3s) holds, then (H6)$\Leftrightarrow$$L_\cdot^+\in \s^p$;
\item [(iv)] (H1)+(H2w)+$g(\cdot,0,0)\in\hcal$ $\Rightarrow$(A); (H1)+(HH)$\Rightarrow$(A).
\end{itemize}

We only show (iv). In fact, if $g$ satisfies (H1) and (H2w) with $g(\cdot,0,0)\in\hcal$, then in view of Remark 1, it follows that $\as$, for each $y\in\R$, $z\in\R^d$,\vspace{-0.2cm}
$$
\begin{array}{lll}
\Dis g(\cdot,y,z){\rm sgn}(y)&\leq &\Dis |(g(\cdot,y,z)-g(\cdot,y,0)){\rm sgn}(y)|\\
&& \Dis +(g(\cdot,y,0)-g(\cdot,0,0)){\rm sgn}(y)+|g(\cdot,0,0)|\\
&\leq &\Dis f_\cdot+\mu |y|+\lambda |z|+\rho(|y|)+|g(\cdot,0,0)|\\
&\leq &f_\cdot+|g(\cdot,0,0)|+A
+(\mu+A)|y|+\lambda |z|.
\end{array}
$$
Hence, $g$ satisfies the assumption (A) with
$$\bar f_\cdot=f_\cdot+|g(\cdot,0,0)|+A,\ \bar\mu=\mu+A\ \ {\rm and}\ \ \bar\lambda=\lambda.$$
Furthermore, if $g$ satisfies assumptions (H1) and (HH), then in view of Remark 1, it follows that $\as$, for each $y\in\R$, $z\in\R^d$,\vspace{-0.2cm}
$$
\begin{array}{lll}
\Dis g(\cdot,y,z){\rm sgn}(y)&\leq &\Dis (g(\cdot,y,z)-g(\cdot,0,z)){\rm sgn}(y)+ |g(\cdot,0,z)|\\
&\leq &\Dis \rho(|y|)+f_\cdot+\lambda |z|\\
&\leq &\Dis f_\cdot+A+A|y|+\lambda |z|.
\end{array}
$$
Hence, $g$ satisfies the assumption (A) with
$$\bar f_\cdot=f_\cdot+A,\ \bar\mu=A\ \ {\rm and}\ \ \bar\lambda=\lambda.$$

\section{A priori estimates}

By virtue of It\^{o}'s formula, the Burkholder-Davis-Gundy (BDG for short) inequality and H\"{o}lder's inequality as well as the stopping time technique and Fatou's Lemma, using a similar argument as that in the proof of Proposition 2.4 of \citet{Izu13} we can prove the following lemma 1. The proof is omitted here.\vspace{0.2cm}

{\bf Lemma 1}\ Let $(\bar Y_\cdot,\bar Z_\cdot,\bar V_\cdot)\in \s\times\M\times\vcal$ satisfy the following equation:
\begin{equation}
\bar Y_t=\bar Y_T+\int_t^T {\rm d}\bar V_s-\int_t^T \bar Z_s{\rm d}B_s,\ \ t\in \T.
\end{equation}
We have
\begin{itemize}
\item [(i)] If $\bar Y_\cdot\in \s^p$ for some $p>0$, then there exists a constant $C_1>0$ depending only on $p$ such that for each $t\in\T$ and $(\F_t)$-stopping time $\tau$ valued in $\T$,
$$
\begin{array}{ll}
&\Dis\E\left[\left.\left(\int_{t\wedge\tau}
^{T\wedge\tau}|\bar Z_s|^2{\rm d}s\right)^{p\over 2}\right|\F_t\right]\\
\leq & \Dis C_1\E\left[\left.\sup\limits_{s\in [t,T]}|\bar Y_{s\wedge\tau}|^p+\sup\limits_{s\in [t,T]}\left[\left(\int_{s\wedge\tau}^{T\wedge\tau} \bar Y_r{\rm d}\bar V_r\right)^+\right]^{p\over 2}\right|\F_t\right].
\end{array}
$$

\item [(ii)] If $\bar Y_\cdot\in \s^p$ for some $p>1$, then there exists a constant $C_2>0$ depending only on $p$ such that for each $t\in\T$ and $(\F_t)$-stopping time $\tau$ valued in $\T$,
$$
\begin{array}{ll}
&\Dis\E\left[\left.\sup\limits_{s\in [t,T]}|\bar Y_{s\wedge\tau}|^p+
\int_{t\wedge\tau}^{T\wedge\tau} |\bar Y_s|^{p-2}\mathbbm{1}_{\{|\bar Y_s|\neq 0\}}|\bar Z_s|^2{\rm d}s
\right|\F_t\right]\vspace{0.1cm}\\
\leq &\Dis C_2\E\left[\left.|\bar Y_\tau|^p+\sup\limits_{s\in [t,T]}\left(\int_{s\wedge\tau}^{T\wedge\tau} |\bar Y_r|^{p-1}{\rm sgn}(\bar Y_r){\rm d}\bar V_r\right)^+\right|\F_t\right].
\end{array}
$$
\end{itemize}

By virtue of Lemma 1 we can prove the following Lemma 2. The proof is classical, see, for example, the proof of Lemma 3.1 and Proposition 3.2 in \citet{Bri03}, we omit it here.

{\bf Lemma 2}\ Assume that the assumption (A) is satisfied for the generator $g$. Let $(Y_\cdot,Z_\cdot,V_\cdot)\in \s\times\M\times\vcal$ satisfy the following equation:
$$
Y_t=Y_T+\int_t^T g(s,Y_s,Z_s){\rm d}s+\int_t^T {\rm d}V_s-\int_t^T Z_s{\rm d}B_s,\ \ t\in \T,
$$
If $Y_\cdot\in \s^p$ and $(V_\cdot,\bar f_\cdot)\in \vcal^p\times\hcal^p$ for some $p>1$, then $Z_\cdot\in \M^p$, and there exists a constant $C>0$ depending only on $p,\bar\mu,\bar\lambda,T$ such that for each $t\in\T$,
$$
\Dis \E\left[\left.\sup\limits_{s\in [t,T]}|Y_s|^p+\left(\int_t^T|Z_s|^2{\rm d}s\right)^{p\over 2}\right|\F_t\right]\\
\leq \Dis C\E\left[\left.|Y_T|^p+|V|^p_T+\left(\int_t^T \bar f_s\ {\rm d}s\right)^p
\right|\F_t\right].\vspace{0.2cm}
$$

{\bf Remark 3}\ Note that in case of $t=0$, Lemma 2 has been obtained in Proposition 3.5 of \citet{Kli13}. \vspace{0.2cm}

By Lemma 1 we can also deduce the following important a priori estimate. \vspace{0.2cm}

{\bf Lemma 3}\ Let $(Y_\cdot,Z_\cdot,V_\cdot,K_\cdot)\in \s\times\M\times\vcal\times\vcal^+$ satisfy the following equation:
$$
Y_t=Y_T+\int_t^T g(s,Y_s,Z_s){\rm d}s+\int_t^T {\rm d}V_s+\int_t^T {\rm d}K_s-\int_t^T Z_s{\rm d}B_s,\ \ t\in \T,
$$
and let $p>0$. We have
\begin{itemize}
\item [(i)] Assume that the assumption (A) is satisfied for the generator $g$. If $Y_\cdot\in \s^p$, then there exists a nonnegative constant $C$ depending only on $p,\bar\mu,\bar\lambda,T$ such that for each $t\in\T$,
\begin{equation}
\begin{array}{ll}
&\Dis \E\left[\left.\left(\int_t^T|Z_s|^2{\rm d}s\right)^{p\over 2}\right|\F_t\right]\\
\leq &\Dis C\E\left[\left.\sup\limits_{s\in [t,T]}|Y_s|^p+|V|^p_T+\left(\int_t^T \bar f_s\ {\rm d}s\right)^p+\left(\int_t^T |Y_s|{\rm d}K_s\right)^{p\over 2}
\right|\F_t\right]
\end{array}
\end{equation}
and
\begin{equation}
\begin{array}{ll}
&\Dis \E\left[\left.\left(\int_t^T |g(s,Y_s,Z_s)|{\rm d}s\right)^p\right|\F_t\right]\vspace{0.1cm}\\
\leq &\Dis C\E\left[\sup\limits_{s\in [t,T]}|Y_s|^p+|V|^p_T+\left(\int_t^T \bar f_s\ {\rm d}s\right)^p+|K_T-K_t|^p
\right.\\
&\hspace{1.0cm} \Dis\left.\left.+\left(\int_t^T|Z_s|^2{\rm d}s\right)^{p\over 2}
\right|\F_t\right].
\end{array}
\end{equation}

\item [(ii)] Assume that the following assumption (B) holds:
\begin{itemize}
\item [{\bf (B)}] There exists two nonnegative constants $\tilde\mu$ and $\tilde\lambda$ such that $\as$,
    $$
     g(t,Y_t,Z_t)\geq -\left(\tilde f_t+\tilde\mu|Y_t|+\tilde\lambda|Z_t|\right),
    $$
    where $\tilde f_t$ is a nonnegative process belonging to $\hcal$.
\end{itemize}
Then there exists a nonnegative constant $\tilde C$ depending only on $p,\tilde\mu,\tilde\lambda,T$ such that for each $t\in\T$,
\begin{equation}
\begin{array}{ll}
&\Dis \E\left[\left.|K_T-K_t|^p\right|\F_t\right]\\
\leq &\Dis \tilde C\E\left[\left.\sup\limits_{s\in [t,T]}|Y_s|^p+|V|^p_T+\left(\int_t^T \tilde f_s\ {\rm d}s\right)^p+\left(\int_t^T|Z_s|^2{\rm d}s\right)^{p\over 2}
\right|\F_t\right].
\end{array}
\end{equation}

\item [(iii)] Assume that $g$ satisfies (A), and that (B) holds. If $Y_\cdot\in \s^p$ and
\begin{equation}
(V_\cdot,\bar f_\cdot,\tilde f_\cdot)\in \vcal^p\times\hcal^p\times\hcal^p,
\end{equation}
then $(Z_\cdot,K_\cdot)\in \M^p\times\vcal^{+,p}$, and there exists a nonnegative constant $\bar C$ depending only on $p,\bar\mu,\bar\lambda, \tilde\mu,\tilde\lambda,T$ such that for each $t\in\T$,
\begin{equation}
\begin{array}{ll}
&\Dis \E\left[\left.\left(\int_t^T|Z_s|^2{\rm d}s\right)^{p\over 2}+|K_T-K_t|^p+\left(\int_t^T|g(s,Y_s,Z_s)|{\rm d}s\right)^p\right|\F_t\right]\\
\leq &\Dis \bar C\E\left[\left.\sup\limits_{s\in [t,T]}|Y_s|^p+|V|^p_T+\left(\int_t^T \bar f_s\ {\rm d}s\right)^p+\left(\int_t^T \tilde f_s\ {\rm d}s\right)^p
\right|\F_t\right].
\end{array}
\end{equation}
\end{itemize}

{\bf Proof.}\ (i) Observe that
$$(\bar Y_\cdot,\bar Z_\cdot,\bar V_\cdot):=\left(Y_\cdot,Z_\cdot,\int_0^\cdot g(s,Y_s,Z_s){\rm d}s+V_\cdot+K_\cdot\right)$$
satisfies equation (2). It follows from (i) of Lemma 1 that if $\bar Y_\cdot\in \s^p$, then there exists a constant $C_1>0$ depending only on $p$ such that for each $t\in\T$,
\begin{equation}
\E\left[\left.\left(\int_t^T|Z_s|^2{\rm d}s\right)^{p\over 2}\right|\F_t\right]\leq C_1\E\left[\left.\sup\limits_{s\in [t,T]}|Y_s|^p+\sup\limits_{s\in [t,T]}\left[\left(\int_s^T Y_r{\rm d}\bar V_r\right)^+\right]^{p\over 2}\right|\F_t\right].
\end{equation}
Furthermore, it follows from (A) and H\"{o}lder inequality that for each $0\leq t\leq s\leq T$,
$$
\begin{array}{lll}
\Dis \left(\int_s^T Y_r{\rm d}\bar V_r\right)^+&\leq & \Dis \bar\mu \sup_{r\in [t,T]}|Y_r|^2+\sup_{r\in [t,T]}|Y_r|\left[\int_t^T \bar f_r\ {\rm d}r+\bar\lambda T
\left(\int_t^T|Z_r|^2{\rm d}r\right)^{1\over 2}\right]\\
&&\Dis +\sup_{r\in [t,T]}|Y_r||V|_T+\int_t^T |Y_r|{\rm d}K_r,
\end{array}
$$
and then by virtue of the basic inequalities
$$ab\leq (a^2+b^2)/2\ \ {\rm and}\ \ (|a|+|b|)^q\leq 2^q(|a|^q+|b|^q),\ q>0,$$
we can get the existence of a constant $C_2>0$ depending only on $p,\bar \mu,\bar\lambda,T$ such that for each $t\in\T$,
\begin{equation}
\begin{array}{lll}
\Dis C_1\sup\limits_{s\in [t,T]}\left[\left(\int_s^T Y_r{\rm d}\bar V_r\right)^+\right]^{p\over 2}&\leq&\Dis {1\over 2}\left(\int_t^T|Z_s|^2{\rm d}s\right)^{p\over 2}+C_2\sup\limits_{s\in [t,T]}|Y_s|^p+C_2|V|^p_T\\
&&\Dis +C_2\left(\int_t^T \bar f_s\ {\rm d}s\right)^p+C_2\left(\int_t^T |Y_s|{\rm d}K_s\right)^{p\over 2}.
\end{array}
\end{equation}
Thus, if $Z_\cdot\in\M^p$, then (3) follows from (8) and (9). Otherwise, for each positive integer $k\geq 1$, define the following $(\F_t)$-stopping time:
\begin{equation}
\tau_k:=\inf\{t\in\T:\ \ \int_0^t |Z_s|^2{\rm d}s\geq k\}\wedge T.
\end{equation}
Note that $\tau_k\To T$ as $k\To +\infty$ due to the fact that $Z_\cdot\in \M$. In the above argument beginning from (8) till (9), replacing respectively
$$
\int_t^T,\ \int_s^T,\ \sup\limits_{s\in [t,T]}|Y_s|^p,\  \sup\limits_{r\in [t,T]}|Y_r|^2,\  \sup\limits_{r\in [t,T]}|Y_r|,\ |V|_T
$$
with
$$
\int_{t\wedge\tau_k}
^{T\wedge\tau_k},\ \int_{s\wedge\tau_k}
^{T\wedge\tau_k},\ \sup\limits_{s\in [t,T]}|Y_{s\wedge\tau_k}|^p,\  \sup\limits_{r\in [t,T]}|Y_{r\wedge\tau_k}|^2,\  \sup\limits_{r\in [t,T]}|Y_{r\wedge\tau_k}|,\ |V|_{T\wedge\tau_k}
$$
yields that for each $k\geq 1$ and each $t\in\T$,
$$
\begin{array}{ll}
&\Dis \E\left[\left.\left(\int_{t\wedge\tau_k}
^{T\wedge\tau_k}|Z_s|^2{\rm d}s\right)^{p\over 2}\right|\F_t\right]\\
\leq &\Dis C\E\left[\left.\sup\limits_{s\in [t,T]}|Y_{s\wedge\tau_k}|^p+|V|^p_{T\wedge\tau_k}+
\left(\int_{t\wedge\tau_k}
^{T\wedge\tau_k} \bar f_s\ {\rm d}s\right)^p+\left(\int_{t\wedge\tau_k}
^{T\wedge\tau_k} |Y_s|{\rm d}K_s\right)^{p\over 2}\right|\F_t\right],
\end{array}
$$
where the constant $C>0$ depends only on $p,\bar \mu,\bar\lambda,T$. Thus, letting $k\To +\infty$ in the above inequality and using Fatou's lemma we get (3).

In the sequel, it follows from (A) that $\as$,
\begin{equation}
\begin{array}{lll}
\Dis |g(\cdot,Y_\cdot,Z_\cdot)| &=& \Dis |-{\rm sgn}(Y_\cdot)g(\cdot,Y_\cdot,Z_\cdot)|\\
&\leq & |f_\cdot+\bar\mu |Y_\cdot|+\bar\lambda |Z_\cdot|-{\rm sgn}(Y_\cdot)g(\cdot,Y_\cdot,Z_\cdot)|+f_\cdot+\bar\mu |Y_\cdot|+\bar\lambda |Z_\cdot|\\
&=& 2\left(f_\cdot+\bar\mu |Y_\cdot|+\bar\lambda |Z_\cdot|\right)-{\rm sgn}(Y_\cdot)g(\cdot,Y_\cdot,Z_\cdot).
\end{array}
\end{equation}
On the other hand, by It\^{o}-Tanaka's formula we know that for each $t\in\T$,
\begin{equation}
\begin{array}{ll}
&\Dis -\int_t^T{\rm sgn}(Y_s)g(s,Y_s,Z_s){\rm d}s\\
\leq &\Dis |Y_T|-|Y_t|+\int_t^T {\rm sgn}(Y_s){\rm d}V_s+\int_t^T {\rm sgn}(Y_s){\rm d}K_s-\int_t^T {\rm sgn}(Y_s)Z_s{\rm d}B_s\\
\leq &\Dis |Y_T|+|V|_T+|K_T-K_t|+\sup_{r\in[t,T]}\left|\int_r^T {\rm sgn}(Y_s)Z_s{\rm d}B_s\right|.
\end{array}
\end{equation}
Thus, (4) follows by combining (11) and (12) and using H\"{o}lder's inequality as well as the BDG inequality.\vspace{0.2cm}

(ii)\ It follows from (B) and H\"{o}lder's inequality that for each $t\in\T$,
$$
\begin{array}{lll}
|K_T-K_t|&=&\Dis Y_t-Y_T-\int_t^Tg(s,Y_s,Z_s){\rm ds}-\int_t^T{\rm d}V_s+\int_t^T Z_s{\rm d}B_s\\
&\leq&\Dis Y_t-Y_T+\int_t^T(\tilde f_s+\tilde\mu |Y_s|+\tilde\lambda |Z_s|){\rm ds}+|V|_T+\left|\int_t^T Z_s{\rm d}B_s\right|\\
&\leq&\Dis (2+\tilde\mu T)\sup\limits_{s\in [t,T]}|Y_s|+|V|_T+\int_t^T \tilde f_s{\rm d}s\\
&&\Dis +\tilde\lambda T\left(\int_t^T|Z_s|^2{\rm d}s\right)^{1\over 2}+\sup_{r\in[t,T]}\left|\int_r^T Z_s{\rm d}B_s\right|,
\end{array}
$$
from which (5) follows immediately by using the basic inequality
$$(\sum\limits_{k=1}^{n}|a_k|)^p\leq n^{p} \sum\limits_{k=1}^{n}|a_k|^p$$
and the BDG inequality.\vspace{0.2cm}

(iii) Since both (A) and (B) are satisfied and $Y_\cdot\in \s^p$, then it follows from (i) and (ii) that both (3), (4) and (5) hold true with constants $C$ and $\tilde C$ respectively. It follows from the basic inequality $2ab\leq \epsilon a^2+b^2/\epsilon,\ \epsilon>0$ with $\epsilon =\tilde C$ that
$$
\begin{array}{lll}
\Dis C\left(\int_t^T |Y_s|{\rm d}K_s  \right)^{p\over 2}&\leq &\Dis C\sup\limits_{s\in\T}|Y_s|^{p\over 2}|K_T-K_t|^{p\over 2}\\
&\leq &\Dis {\tilde C\over 2}C^2\sup\limits_{s\in\T}|Y_s|^p+{1\over 2\tilde C}|K_T-K_t|^p.
\end{array}
$$
Combining (3), (5) and the above inequality yields the existence of a constant $\bar C>0$ depending only on $p,\bar\mu,\bar\lambda, \tilde\mu,\tilde\lambda,T$ such that for each $t\in\T$,
\begin{equation}
\begin{array}{ll}
&\Dis \E\left[\left.\left(\int_t^T|Z_s|^2{\rm d}s\right)^{p\over 2}\right|\F_t\right]\\
\leq &\Dis \bar C\E\left[\left.\sup\limits_{s\in [t,T]}|Y_s|^p+|V|^p_T+\left(\int_t^T \bar f_s\ {\rm d}s\right)^p+\left(\int_t^T \tilde f_s\ {\rm d}s\right)^p
\right|\F_t\right],
\end{array}
\end{equation}
provided that $Z_\cdot\in \M^p$. Using the stopping time technique and Fatou's lemma, similar to the argument beginning from (10) to the end of that paragraph, we can deduce that (13) holds still true for $Z_\cdot$ which only belongs to $\M$. Thus, (7) follows from (5), (4) and (13), and in view of (6), $(Z_\cdot,K_\cdot)\in \M^p\times\vcal^{+,p}$. \vspace{0.2cm}\hfill $\Box$

By virtue of Lemmas 2 and 3 we can prove the following Propositions 1 and 2, which will play important roles later. \vspace{0.1cm}

{\bf Proposition 1}\ Assume that $p>1$, $\xi\in\Lp$ and $V_\cdot\in \vcal^p$. For each $n\geq 1$, suppose that generators $g_n$ satisfy (H1) and (HH) with the same $\rho(\cdot)$, $f_\cdot$, $\psi_\cdot(r)$ and $\lambda$, and let $(Y_\cdot^n,Z_\cdot^n)\in \s^p\times\M^p$ satisfy
$$
Y_t^n=\xi+\int_t^T g_n(s,Y_s^n,Z_s^n){\rm d}s+\int_t^T {\rm d}V_s-\int_t^T Z_s^n{\rm d}B_s,\ \ t\in \T.
$$
Then, there exists a nonnegative constant $C$ depending only on $p,A,\lambda,T$ such that for each $t\in\T$ and $n\geq 1$,
$$
\begin{array}{ll}
&\Dis \E\left[\left.\sup\limits_{s\in [t,T]}|Y_s^n|^p+\left(\int_t^T|Z_s^n|^2{\rm d}s\right)^{p\over 2}+\left(\int_t^T|g_n(s,Y_s^n,Z_s^n)|{\rm d}s\right)^p\right|\F_t\right]\\
\leq &\Dis C\E\left[\left.|\xi|^p+|V|^p_T+\left(\int_t^T f_s{\rm d}s\right)^p+1\right|\F_t\right].
\end{array}
$$

{\bf Proof.} It follows from (iv) of Remark 2 that all $g_n$ satisfy (A) with the same
$$\bar f_\cdot=f_\cdot+A,\ \bar\mu=A\ \ {\rm and}\ \ \bar\lambda=\lambda.$$
Then by Lemma 2 we know that there exists a constant $C_2>0$ depending only on $p,A,\lambda,T$ such that for each $t\in\T$ and $n\geq 1$,
$$
\begin{array}{ll}
&\Dis \E\left[\left.\sup\limits_{s\in [t,T]}|Y_s^n|^p+\left(\int_t^T|Z_s^n|^2{\rm d}s\right)^{p\over 2}\right|\F_t\right]\\
\leq &\Dis C_2\E\left[\left.|\xi|^p+|V|^p_T+\left(\int_t^T f_s{\rm d}s\right)^p+1\right|\F_t\right],
\end{array}
$$
form which together with (4) in Lemma 3 the conclusion follows. \vspace{0.2cm} $\Box$\hfill

{\bf Proposition 2}\ Assume that $p>1$, $\xi\in\Lp$ and $V_\cdot\in \vcal^p$. For $n\geq 1$, suppose that generators $g_n$ and $g$ satisfy (H1) with the same $\rho(\cdot)$, $g(\cdot,0,0)\in \hcal^p$ and that there exist two nonnegative constants $\mu, \lambda$ and a nonnegative process $f_\cdot\in \hcal^p$ such that $\as$, $\RE\ y\in\R, z\in\R^{d}$,
\begin{equation}
\Dis |g_n(\omega,t,y,z)-g(\omega,t,y,0)|\leq f_t(\omega)+\mu |y|+\lambda |z|.\vspace{0.1cm}
\end{equation}
For each $n\geq 1$, let $(Y_\cdot^n,Z_\cdot^n,K_\cdot^n)\in \s^p\times\M^p\times\vcal^{+,p}$ satisfy
$$
Y_t^n=\xi+\int_t^T g_n(s,Y_s^n,Z_s^n){\rm d}s+\int_t^T {\rm d}V_s+\int_t^T {\rm d}K_s^n-\int_t^T Z_s^n{\rm d}B_s,\ \ t\in \T.
$$
If the following assumption (C) holds:
\begin{itemize}
\item [{\bf (C)}] There exists a process $\bar X_\cdot\in \s^p$ such that $g(\cdot,\bar X_\cdot,0)\in \hcal^p$ and $\bar X_t\geq Y_t^n$ for each $t\in \T$ and $n\geq 1$,\vspace{-0.1cm}
\end{itemize}
then there exists a nonnegative constant $C$ depending only on $p,\mu, \lambda,A,T$ such that for each $t\in\T$ and $n\geq 1$,
\begin{equation}
\begin{array}{ll}
&\Dis\E\left[\left.\left(\int_t^T|Z_s^n|^2{\rm d}s\right)^{p\over 2}+|K_T^n-K_t^n|^p+\left(\int_t^T|g_n(s,Y_s^n,Z_s^n)|{\rm d}s\right)^p\right|\F_t\right]\\
\leq &\Dis C\E\left[\sup\limits_{s\in [t,T]}|Y_s^n|^p+|V|^p_T+\sup\limits_{s\in [t,T]}|\bar X_s|^{p}+\left(\int_t^T f_s{\rm d}s\right)^p+1\right.\\
&\hspace{1cm}\Dis+\left.\left.\left(\int_t^T |g(s,\bar X_s,0)|\ {\rm d}s\right)^p+\left(\int_t^T |g(s,0,0)|\ {\rm d}s\right)^p\right|\F_t\right].
\end{array}
\end{equation}

{\bf Proof.}\ In view of Remark 1, it follows from (14) together with (H1) for $g$ that $\as$, for each $y\in\R$, $z\in\R^d$ and $n\geq 1$,\vspace{-0.2cm}
$$
\begin{array}{lll}
\Dis g_n(\cdot,y,z){\rm sgn}(y)&\leq &\Dis |(g_n(\cdot,y,z)-g(\cdot,y,0)){\rm sgn}(y)| \\
&&\Dis +(g(\cdot,y,0)-g(\cdot,0,0)){\rm sgn}(y)+|g(\cdot,0,0))|\\
&\leq &\Dis f_\cdot+\mu |y|+\lambda |z|+\rho(|y|)+|g(\cdot,0,0))|\\
&\leq &\Dis f_\cdot+|g(\cdot,0,0))|+A+(\mu+A)|y|+\lambda |z|.
\end{array}
$$
Hence, for each $n\geq 1$, $g_n$ satisfies the assumption (A) with the same
$$\bar f_\cdot=f_\cdot+|g(\cdot,0,0)|+A,\ \bar\mu=\mu+A\ \ {\rm and}\ \ \bar\lambda=\lambda.$$
Furthermore, note by (C) that $\bar X_t\geq Y_t^n$ for each $t\in \T$ and $n\geq 1$. Once again, in view of Remark 1, it follows from (H1) for $g_n$ together with (14) that $\as$, for each $n\geq 1$,\vspace{-0.2cm}
$$
\begin{array}{lll}
\Dis -g_n(\cdot,Y_\cdot^n,Z_\cdot^n)&\leq &\Dis
\left(g_n(\cdot,\bar X_\cdot,Z_\cdot^n)
-g_n(\cdot,Y_\cdot^n,Z_\cdot^n)\right)\\
&&\Dis+\left|g(\cdot,\bar X_\cdot,0)-g_n(\cdot,\bar X_\cdot,Z_\cdot^n)
\right|+|g(\cdot,\bar X_\cdot,0)|\\
&\leq &\Dis \rho(|\bar X_\cdot-Y_\cdot^n|)+f_\cdot+\mu |\bar X_\cdot|+\lambda |Z_\cdot^n|+|g(\cdot,\bar X_\cdot,0)|\\
&\leq &\Dis |g(\cdot,\bar X_\cdot,0)|+(\mu+A)|\bar X_\cdot|+f_\cdot+A+A|Y_\cdot^n|+\lambda |Z_\cdot^n|.
\end{array}
$$
Hence, the assumption (B) holds also true for each $g_n$ with the same
$$
\tilde f_\cdot=|g(\cdot,\bar X_\cdot,0)|+(\mu+A)|\bar X_\cdot|+f_\cdot+A,\ \tilde\mu=A\ \ {\rm and}\ \ \tilde\lambda=\lambda.
$$
Thus, (15) follows from (iii) of Lemma 3.\vspace{0.2cm} $\Box$ \hfill

\section{Penalization, approximation and comparison theorem}

In this section, we will put forward and prove the convergence of sequence of $L^p$ solutions for penalization and approximation equations of RBSDE (1) with $L^p\ (p>1)$ data under assumption (HH) and some very elementary conditions, and a comparison theorem for $L^p$ solutions of RBSDEs with $L^p\ (p>1)$ data under assumptions (H1) and (H2).

{\bf Proposition 3}\ (Penalization)\ Assume that the generator $g$ satisfies (H4) and (HH) with $f_\cdot, \psi_\cdot(r)$ and $\lambda$. Let $p>1$, $V_\cdot\in \vcal^p$ and (H5) be satisfied for $\xi$ and $L_\cdot$. For $n\geq 1$, assume that $(Y_\cdot^n,Z_\cdot^n)\in \s^p\times\M^p$ satisfies the following penalization BSDE:
\begin{equation}
Y_t^n=\xi+\int_t^T g(s,Y_s^n,Z_s^n){\rm d}s+\int_t^T {\rm d}V_s+\int_t^T {\rm d}K_s^n-\int_t^T Z_s^n{\rm d}B_s,\ \ t\in \T
\end{equation}
with
\begin{equation}
K_t^n:=n\int_0^t\left(Y_s^n-L_s\right)^-{\rm d}s,\ \ t\in\T.\vspace{0.1cm}
\end{equation}
If $Y_\cdot^n$ increases in $n$ and there exists a random variable $\eta\in {\mathbb L}^1(\F_T)$ such that for each $n\geq 1$ and $t\in\T$,
\begin{equation}
\begin{array}{l}
\Dis\E\left[\left.\sup\limits_{s\in [t,T]}|Y_s^n|^p+\left(\int_t^T|Z_s^n|^2{\rm d}s\right)^{p\over 2}+|K_T^n-K_t^n|^p+\left(\int_t^T|g(s,Y_s^n,Z_s^n)|{\rm d}s\right)^p\right|\F_t\right]\\
\ \ \ \ \leq \E\left[\left.\eta\right|\F_t\right],
\end{array}
\end{equation}
then there exists a triple $(Y_\cdot,Z_\cdot,K_\cdot)\in \s^p\times\M^p\times\vcal^{+,p}$ which solves RBSDE (1),
$$\lim\limits_{n\To \infty}\left(\|Y_\cdot^n-Y_\cdot\|_{\s^p}+
\|Z_\cdot^n-Z_\cdot\|_{\M^p}\right)=0,$$
and there exists a subsequence $\{K_\cdot^{n_j}\}$ of $\{K_\cdot^n\}$ such that
$$\lim\limits_{j\To\infty}\sup\limits_{t\in\T}
|K_t^{n_j}-K_t|=0.$$

{\bf Proof.}\ Since $Y_\cdot^n$ increases in $n$, there exists a process $Y_\cdot$ such that $\ps$, $Y_t^n\uparrow Y_t$ for each $t\in\T$. By Fatou's lemma and (18) we can deduce that
\begin{equation}
\begin{array}{lll}
\Dis\E\left[\sup\limits_{t\in\T}|Y_t|^p\right]
&=&\Dis\E\left[\sup\limits_{t\in\T}\liminf
\limits_{n\To\infty}|Y_t^n|^p\right]
\leq \E\left[\liminf\limits_{n\To\infty}
\sup\limits_{t\in\T}|Y_t^n|^p\right]\vspace{0.1cm}\\
&\leq& \Dis \liminf\limits_{n\To\infty}\E\left[\sup\limits_{t\in\T}|Y_t^n|^p\right]\leq \E\left[\eta\right]<+\infty.
\end{array}
\end{equation}
Furthermore, by (18) we can also get that
\begin{equation}
\sup\limits_{n\geq 1}|Y^n_t|\leq \left(\E\left[\left.\eta\right|\F_t\right]
\right)^{1\over p},\ \ t\in\T
\end{equation}
and
\begin{equation}
\sup\limits_{n\geq 1}\E\left[\left(\int_0^T|Z_t^n|^2{\rm d}t\right)^{p\over 2}+|K_T^n|^p+\left(\int_0^T|g(t,Y_t^n,Z_t^n)|{\rm d}t\right)^p\right]\leq \E\left[\eta\right]<+\infty.\vspace{0.1cm}
\end{equation}
The rest proof is divided into 6 steps.\vspace{0.2cm}

{\bf Step 1.}\ We show that $Y_\cdot$ is a c\`{a}dl\`{a}g process. For each integer $l,q\geq 1$, introduce the following two $(\F_t)$-stopping times:
$$
\begin{array}{rll}
\tau_l&:=&\Dis \inf\left\{t\geq 0:\  \left(\E\left[\left.\eta\right|\F_t\right]
\right)^{1\over p}+\int_0^t f_s{\rm d}s+L_t\geq l\right\}\wedge T;\\
\sigma_{l,q}&:=&\Dis \inf\left\{t\geq 0:\  \int_0^t \psi_s(l) {\rm d}s\geq q\right\}\wedge \tau_l.
\end{array}
$$
Then we have, $\tau_l\To T$ as $l\To \infty$, $\sigma_{l,q}\To \tau_l$ as $q\To \infty$ for each $l\geq 1$,
$$
\mathbb{P}\left(\left\{\omega:\ \exists l_0(\omega)\geq 1, \ \RE l\geq l_0(\omega), \ \tau_l(\omega)=T\right\}\right)=1
$$
and
\begin{equation}
\mathbb{P}\left(\left\{\omega:\ \exists l_0(\omega), q_0(\omega)\geq 1,\ \RE l\geq l_0(\omega), \RE q\geq q_0(\omega), \ \sigma_{l,q}(\omega)=T\right\}\right)=1.
\end{equation}

Now, let us arbitrarily fix a pair of $l,q\geq 1$. Since $g$ satisfies (HH) with $f_\cdot, \psi_\cdot(r)$ and $\lambda$, and (20) is satisfied, in view of the definitions of $\tau_l$ and $\sigma_{l,q}$, we know that $\as$, for each $n\geq 1$,
\begin{equation}
|h^{n;l,q}_\cdot|\leq \mathbbm{1}_{\cdot\leq \tau_{l}}f_\cdot+\mathbbm{1}_{\cdot\leq \sigma_{l,q}}\psi_\cdot(l)+\lambda |Z_\cdot^n|\vspace{0.1cm}
\end{equation}
with $h^{n;l,q}_\cdot:=\mathbbm{1}_{\cdot\leq \sigma_{l,q}}g(\cdot,Y_\cdot^n,Z_\cdot^n)$,
\begin{equation}
\E\left[\int_0^T \mathbbm{1}_{t\leq \tau_{l}}f_t{\rm d}t\right]\leq l\ \ {\rm and}\ \ \E\left[\int_0^T \mathbbm{1}_{t\leq \sigma_{l,q}}\psi_t(l){\rm d}t\right]\leq q,
\end{equation}
from which together with (21), we can deduce that there exists a subsequence $\{h^{n_j;l,q}_\cdot\}_{j=1}^{\infty}$ of the sequence $\{h^{n; l,q}_\cdot\}_{n=1}^{\infty}$ which converges weakly to a process $h^{l,q}_\cdot$ in $\hcal^1$.  Now, take any bounded
linear functional $\Phi(\cdot)$ defined on $\mathbb{L}^1(\F_T)$. Then there exists a constant $b>0$ such that for each $\overline{h}_\cdot \in \hcal^1$ and every $(\F_t)$-stopping time $\bar\tau$ valued in $\T$, we have
$$
\left|\Phi(\int_0^{\bar\tau}\overline h_s {\rm d}s)\right|\leq b\left\|\int_0^{\bar\tau}\overline h_s {\rm d}s\right\|_{\mathbb{L}^1}\leq b\left\|\overline h\right\|_{\hcal^1}.
$$
Hence, for each $(\F_t)$-stopping time $\bar\tau$ valued in $\T$, $\Phi(\int_0^{\bar\tau}\cdot \ {\rm d}s )$ is a bounded linear functional defined on $\hcal^1$, which means that
$$\lim\limits_{j\To \infty} \Phi\left(\int_0^{\bar\tau} h_s^{n_j;l,q} {\rm d}s \right)=\Phi\left(\int_0^{\bar\tau} h^{l,q}_s {\rm d}s \right).$$
As a result, for every $(\F_t)$-stopping time $\tau$ with $0\leq\tau\leq \sigma_{l,q}$, as $j\To \infty$,
\begin{equation}
\int_0^\tau g(s,Y_s^{n_j},Z_s^{n_j}) {\rm d}s=\int_0^\tau h_s^{n_j;l,q} {\rm d}s\ \To\  \int_0^\tau h^{l,q}_s {\rm d}s\ \ {\rm weakly\ in}\ \mathbb{L}^1(\F_T).
\end{equation}
Furthermore, it follows from (21) and Lemma 4.4 of \citet{Kli12} that there exists a process $Z_\cdot\in \M^p$ and a subsequence of the sequence $\{n_j\}_{j=1}^{\infty}$, still denoted by itself, such that for every $(\F_t)$-stopping time $\bar\tau$ valued in $\T$, as $j\To \infty$,
\begin{equation}
\int_0^{\bar\tau} Z_s^{n_j}{\rm d}B_s\To \int_0^{\bar\tau} Z_s{\rm d}B_s\ \ {\rm weakly\ in}\ \mathbb{L}^p(\F_T)\ {\rm and \ then\ in}\ \mathbb{L}^1(\F_T).
\end{equation}
In the sequel, we define
$$
K^{l,q}_t:=Y_0-Y_t-\int_0^t h^{l,q}_s{\rm d}s-\int_0^t {\rm d}V_s+\int_0^t Z_s{\rm d}B_s, \ \ t\in\T.
$$
Then, in view of (25), (26) and the fact that for each $(\F_t)$-stopping time $\bar\tau$ valued in $\T$, $Y^n_{\bar\tau}\uparrow Y_{\bar\tau}$ in $\mathbb{L}^1(\F_T)$, we can get that for every $(\F_t)$-stopping time $\tau$ such that $0\leq\tau\leq \sigma_{l,q}$, the sequence of random variables
$$
K_\tau^{n_j}=Y_0^{n_j}-Y_\tau^{n_j}-\int_0^\tau
g(s,Y_s^{n_j},Z_s^{n_j}){\rm d}s-\int_0^\tau{\rm d}V_s+\int_0^\tau Z_s^{n_j}{\rm d}B_s
$$
converges weakly to $K^{l,q}_\tau$ in $\mathbb{L}^1(\F_T)$ as $j\To \infty$. Consequently, since $K_\cdot^n\in \vcal^+$ for each $n\geq 1$, we know that $$K^{l,q}_{\tau_1\wedge \sigma_{l,q}}\leq K^{l,q}_{\tau_2\wedge \sigma_{l,q}}\vspace{0.2cm}$$
for any $(\F_t)$-stopping times $\tau_1\leq \tau_2$ valued in $\T$, and in view of the definition of $K_\cdot^{l,q}$ and the facts that $Y^n_\cdot\uparrow Y_\cdot$ and $Y^n_\cdot\in \s^p$ for each $n\geq 1$, it is not hard to check that $K^{l,q}_{\cdot}$ is a $(\F_t)$-optional process with $\ps$ upper semi-continuous paths. Thus, Lemma A.3 in \citet{Bay14} yields that $K^{l,q}_{\cdot\wedge \sigma_{l,q}}$ is a nondecreasing process, and then it has $\ps$ right lower semi-continuous paths. Hence, $K^{l,q}_{\cdot\wedge \sigma_{l,q}}$ is c\`{a}dl\`{a}g and so is $Y_{\cdot\wedge \sigma_{l,q}}$ from the definition of $K^{l,q}_\cdot$. Finally, it follows from (22) that $Y_\cdot$ is also a c\`{a}dl\`{a}g process.\vspace{0.2cm}

{\bf Step 2.}\ We show that $Y_t\geq L_t$ for each $t\in \T$ and as $n\To \infty$,
\begin{equation}
\sup\limits_{t\in\T}(Y^n_t-L_t)^-\To 0.
\end{equation}
In fact, it follows from (21) and the definition of $K^n_\cdot$ that for each $n\geq 1$,
$$\E\left[\left(\int_0^T(Y^n_t-L_t)^-{\rm d}t\right)^p\right]\leq {\E[\eta]\over n^p}.$$
Hence, by Fatou's lemma and H\"{o}lder's inequality,
$$
\E\left[\int_0^T(Y_t-L_t)^-{\rm d}t\right]\leq \liminf\limits_{n\To\infty}\E\left[\int_0^T(Y^n_t-L_t)^-{\rm d}t\right]\leq \lim\limits_{n\To\infty}{(\E[\eta])^{1\over p}\over n}=0,
$$
which implies that $$\E\left[\int_0^T(Y_t-L_t)^-{\rm d}t\right]=0.$$ Since $Y_\cdot-L_\cdot$ is a c\`{a}dl\`{a}g process, $(Y_t-L_t)^-=0$ and hence
$Y_t\geq L_t$ for each $t\in [0,T)$. Moreover, $Y_T=Y^n_T=\xi\geq L_T$. Hence $$(Y^n_t-L_t)^-\downarrow 0$$
for each $t\in [0,T]$ and by Dini's theorem, (27) follows. \vspace{0.2cm}

{\bf Step 3.}\ We show the convergence of the sequence $\{Y_\cdot^n\}$. Let $\tau_l$ and $\sigma_{l,q}$ be the sequences of $(\F_t)$-stopping times defined in Step 1. For each $n,m\geq 1$, observe that
\begin{equation}
\begin{array}{lll}
\Dis (\bar Y_\cdot,\bar Z_\cdot,\bar V_\cdot)&
:=&\Dis (Y_\cdot^n-Y_\cdot^m,Z_\cdot^n-Z_\cdot^m,\\
&&\Dis \ \ \int_0^\cdot \left(g(s,Y_s^n,Z_s^n)-g(s,Y_s^m,Z_s^m)\right){\rm d}s+\left(K_\cdot^n-K_\cdot^m\right) )
\end{array}
\end{equation}
satisfies equation (2). It then follows from (ii) of Lemma 1 with $p=2$, $t=0$ and $\tau=\sigma_{l,q}$ that there exists a constant $C>0$ such that for each $n,m,l,q\geq 1$,
\begin{equation}
\begin{array}{ll}
&\Dis\E\left[\sup\limits_{t\in [0,T]} |Y_{t\wedge\sigma_{l,q}}^n
-Y_{t\wedge\sigma_{l,q}}^m|^2\right]\\
\leq &\Dis C\E\left[|Y_{\sigma_{l,q}}^n-Y_{\sigma_{l,q}}^m|^2 +\sup\limits_{t\in [0,T]}\left(\int_{t\wedge \sigma_{l,q}}^{\sigma_{l,q}}(Y^n_s-Y_s^m)
\left({\rm d}K_s^n-{\rm d}K_s^m\right) \right)^+\right.\vspace{0.1cm}\\
&\hspace{1cm}\Dis +\left. \int_{0}^{\sigma_{l,q}}|Y^n_t-Y_t^m|
\left|g(t,Y_t^n,Z_t^n)-g(t,Y_t^m,Z_t^m)\right| {\rm d}t\right].
\end{array}
\end{equation}
Furthermore, by virtue of the definition of $K_\cdot^n$ we know that for each $t\in \T$,
\begin{equation}
\begin{array}{ll}
&\Dis \int_{t\wedge \sigma_{l,q}}^{\sigma_{l,q}}(Y^n_s-Y_s^m)
\left({\rm d}K_s^n-{\rm d}K_s^m\right)\\
= &\Dis \int_{t\wedge \sigma_{l,q}}^{\sigma_{l,q}}\left[(Y^n_s-L_s)
-(Y_s^m-L_s)\right]
{\rm d}K_s^n-\int_{t\wedge \sigma_{l,q}}^{\sigma_{l,q}}\left[(Y^n_s-L_s)
-(Y_s^m-L_s)\right]
{\rm d}K_s^m\\
\leq &\Dis \int_{t\wedge \sigma_{l,q}}^{\sigma_{l,q}}(Y_s^m-L_s)^-
{\rm d}K_s^n+\int_{t\wedge \sigma_{l,q}}^{\sigma_{l,q}}(Y_s^n-L_s)^-{\rm d}K_s^m\\
\leq &\Dis \sup\limits_{t\in\T}(Y_{t\wedge\sigma_{l,q} }^m-L_{t\wedge\sigma_{l,q}})^-|K_T^n|+
\sup\limits_{t\in\T}(Y_{t\wedge\sigma_{l,q} }^n-L_{t\wedge\sigma_{l,q}})^-|K_T^m|.
\end{array}
\end{equation}
Combining (23), (29) and (30) together with H\"{o}lder's inequality yields that
\begin{equation}
\begin{array}{ll}
&\Dis\E\left[\sup\limits_{t\in [0,T]} |Y_{t\wedge\sigma_{l,q}}^n
-Y_{t\wedge\sigma_{l,q}}^m|^2\right]\\
\leq &\Dis C\E\left[|Y_{\sigma_{l,q}}^n-Y_{\sigma_{l,q}}^m|^2 +2\int_0^T |Y^n_t-Y_t^m|\left(\mathbbm{1}_{t\leq \tau_l}f_t+\mathbbm{1}_{t\leq\sigma_{l,q}}
\psi_t(l)\right){\rm d}t\right]\\
&\Dis+C\left(\E\left[ \sup\limits_{t\in\T}\left|(Y_{t\wedge\sigma_{l,q} }^m-L_{t\wedge\sigma_{l,q}})^-\right|^{p\over p-1}\right]\right)^{p-1\over p}\left(\E\left[|K_T^n|^p\right]\right)^{1\over p}\vspace{0.1cm}\\
&\Dis+C\left(\E\left[ \sup\limits_{t\in\T}\left|(Y_{t\wedge\sigma_{l,q} }^n-L_{t\wedge\sigma_{l,q}})^-\right|^{p\over p-1}\right]\right)^{p-1\over p}\left(\E\left[|K_T^m|^p\right]\right)^{1\over p}\\
&\Dis+2C\lambda\left(\E\left[\left(\int_{0}^{
\sigma_{l,q}}|Y^n_t-Y_t^m|^2{\rm d}t\right)^{p\over 2(p-1)}\right]\right)^{p-1\over p}\\
&\Dis \hspace{1cm}\times \left(\E\left[\left(\int_{0}^T\left(|Z_t^n|+
|Z_t^m|\right)^2{\rm d}t\right)^{p\over 2}\right]\right)^{1\over p}.\vspace{0.2cm}
\end{array}
\end{equation}
Thus, note that $Y_t^n\uparrow Y_t$ for each $t\in\T$. In view of the definitions of $\tau_l$ and $\sigma_{l,q}$, (20), (21), (24) and (27), by (31) and Lebesgue's dominated convergence theorem we can deduce that for each $l,q\geq 1$, as $n,m\To\infty$,
$$\E\left[\sup\limits_{t\in [0,T]} |Y_{t\wedge\sigma_{l,q}}^n
-Y_{t\wedge\sigma_{l,q}}^m|^2\right]\To 0,$$
which implies that for each $l,q\geq 1$, as $n,m\To\infty$,
$$\sup\limits_{t\in [0,T]} |Y_{t\wedge\sigma_{l,q}}^n
-Y_{t\wedge\sigma_{l,q}}^m|\To 0\ {\rm in\ probability}\ \mathbb{P}.$$
And, by (22) and the fact that $Y_\cdot^n$ increases in $n$ we know that $\ps$,
\begin{equation}
\sup\limits_{t\in [0,T]} |Y_t^n
-Y_t|\To 0,\ \ {\rm as}\ n\To\infty.
\end{equation}
So, $Y_\cdot$ is a continuous process. Finally, note that $|Y_\cdot^n|\leq |Y_\cdot^1|+|Y_\cdot|$ for each $n\geq 1$ and that (19) is satisfied. From (32) and Lebesgue's dominated convergence theorem it follows that
\begin{equation}
\lim\limits_{n\To\infty}\|Y_\cdot^n-Y_\cdot\|_{\s^p}^p =\lim\limits_{n\To\infty}\E\left[\sup\limits_{t\in [0,T]} |Y_t^n-Y_t|^p\right]=0.
\end{equation}

{\bf Step 4.}\ We show the convergence of the sequence $\{Z_\cdot^n\}$. Note that (28) solves (2). It follows from (i) of Lemma 1 with $t=0$ and $\tau=T$ that there exists a constant $C'>0$ such that for each $m,n\geq 1$,
$$
\begin{array}{ll}
&\Dis\E\left[\left(\int_0^T|Z_t^n-Z_t^m|^2{\rm d}t\right)^{p\over 2}\right]\\
\leq & \Dis C'\E\left[\sup\limits_{t\in [0,T]}|Y_t^n-Y_t^m|^p+\sup\limits_{t\in [0,T]}\left[\left(\int_t^T (Y_s^n-Y_s^m)\left({\rm d}K_s^n-{\rm d}K_s^m\right)\right)^+\right]^{p\over 2}\right]\\
&\Dis + C'\E\left[\left(\int_{0}^T|Y^n_t-Y_t^m|
\left|g(t,Y_t^n,Z_t^n)-g(t,Y_t^m,Z_t^m)\right| {\rm d}t\right)^{p\over 2}\right].
\end{array}
$$
Then, it follows from H\"{o}lder's inequality that for each $m,n\geq 1$,
$$
\begin{array}{ll}
&\Dis\E\left[\left(\int_0^T|Z_t^n-Z_t^m|^2{\rm d}t\right)^{p\over 2}\right]\\
\leq & \Dis C'\E\left[\sup\limits_{t\in [0,T]}|Y_t^n-Y_t^m|^p\right]+\left.C'\left(\E\left[\sup\limits_{t\in [0,T]}|Y_t^n-Y_t^m|^p\right]\right)^{1\over 2}\right\{\left(\E\left[|K_T^n|^p\right]\right)^{1\over 2}\\
&\Dis+\left.\left(\E\left[|K_T^m|^p\right]\right)^{1\over 2}+\left(\E\left[\left(\int_{0}^T
\left(|g(t,Y_t^n,Z_t^n)|+|g(t,Y_t^m,Z_t^m)|\right) {\rm d}t\right)^{p}\right]\right)^{1\over 2}\right\},
\end{array}
$$
from which together with (21), (33) and (26) yields that
\begin{equation}
\lim\limits_{n\To\infty}\|Z_\cdot^n-Z_\cdot\|_{\M^p}^p
=\lim\limits_{n\To\infty}\E\left[\left(\int_0^T|Z_t^n
-Z_t|^2{\rm d}t\right)^{p\over 2}\right]=0.
\end{equation}

{\bf Step 5.}\ We show the convergence of the sequence $\{K_\cdot^n\}$. Let $\tau_l$ and $\sigma_{l,q}$ be the sequences of $(\F_t)$-stopping times defined in Step 1. Since $g$ satisfies (H4), by (23), (24), (21), (32) and (34) we can deduce that there exists a subsequence $\{n_j\}$ of $\{n\}$ such that for each $l,q\geq 1$,
$$\lim\limits_{j\To \infty}\int_0^{\sigma_{l,q}}|g(t,Y_t^{n_j},Z_t^{n_j})-
g(t,Y_t,Z_t)|{\rm d}t=0.$$
Then, in view of (22), we have
\begin{equation}
\lim\limits_{j\To \infty}\sup\limits_{t\in\T}\left|\int_0^t
g(t,Y_t^{n_j},Z_t^{n_j}){\rm d}t-\int_0^t
g(t,Y_t,Z_t){\rm d}t\right|=0.
\end{equation}
Combining (32), (34) and (35) yields that $\ps$, for each $t\in\T$,
$$
K_t^{n_j}=Y_0^{n_j}-Y_t^{n_j}-\int_0^t
g(s,Y_s^{n_j},Z_s^{n_j}){\rm d}s-\int_0^t{\rm d}V_s+\int_0^tZ_s^{n_j}{\rm d}B_s
$$
tends to
$$
K_t:=Y_0-Y_t-\int_0^t
g(s,Y_s,Z_s){\rm d}s-\int_0^t{\rm d}V_s+\int_0^tZ_s{\rm d}B_s
$$
as $j\To \infty$ and that
\begin{equation}
\lim\limits_{j\To\infty}\sup\limits_{t\in\T}
|K_t^{n_j}-K_t|=0.
\end{equation}
Hence, $K_\cdot$ is a continuous process.\vspace{0.2cm}

{\bf Step 6.}\ We show that $K_\cdot\in\vcal^{+,p}$ and $(Y_\cdot,Z_\cdot,
K_\cdot)\in \s^p\times\M^p\times\vcal^{+,p}$ is a solution of RBSDE (1). In fact, by Fatou's lemma with (36) and (21) we get that
$$
\begin{array}{lll}
\Dis\E\left[\sup\limits_{t\in\T}|K_t|^p\right]
&=&\E\left[\lim\limits_{j\To\infty}
\sup\limits_{t\in\T}|K_t^{n_j}|^p\right]\leq \Dis \liminf\limits_{j\To\infty}
\E\left[\sup\limits_{t\in\T}|K_t^{n_j}|^p\right]
\vspace{0.1cm}\\
&\leq& \Dis \sup\limits_{j\geq 1}\E\left[|K_T^{n_j}|^p\right]\leq  \E\left[\eta\right]<+\infty.
\end{array}
$$
Hence, $K_\cdot\in\vcal^{+,p}$ and $(Y_\cdot,Z_\cdot,K_\cdot)\in \s^p\times\M^p\times\vcal^{+,p}$ solves
$$
Y_t=\xi+\int_t^Tg(s,Y_s,Z_s){\rm d}s+\int_t^T{\rm d}V_s+\int_t^T{\rm d}K_s-\int_t^TZ_s {\rm d}B_s,\ \ t\in\T.
$$
By Step 2 we know that $Y_t\geq L_t$ for each $t\in\T$, and then
$$
\int_0^T(Y_t-L_t){\rm d}K_t\geq 0.
$$
On the other hand, in view of (32) and (36), it follows from the definition of $K_\cdot^n$ that
$$
\int_0^T(Y_t-L_t){\rm d}K_t=\lim\limits_{j\To \infty}\int_0^T(Y_t^{n_j}-L_t){\rm d}K_t^{n_j}\leq 0.
$$
Consequently, we have
$$
\int_0^T(Y_t-L_t){\rm d}K_t=0.
$$
Thus, $(Y_\cdot,Z_\cdot,
K_\cdot)\in \s^p\times\M^p\times\vcal^{+,p}$ solves (1.1). Proposition 3 is then proved.\vspace{0.2cm}\hfill $\Box$

{\bf Proposition 4}\ (Approximation)\ Assume that for each $n\geq 1$, the generator $g_n$ satisfies (HH) with the same $f_\cdot, \psi_\cdot(r)$ and $\lambda$. Let $p>1$, $V_\cdot\in \vcal^p$ and (H5) be satisfied for $\xi$ and $L_\cdot$. For $n\geq 1$, assume that $(Y_\cdot^n,Z_\cdot^n,K_\cdot^n)\in \s^p\times\M^p\times\vcal^{+,p}$ is a solution of RBSDE $(\xi,g_n+{\rm d}V,L)$. If $Y_\cdot^n$ increases or decreases in $n$, $g_n$ tends locally uniformly in $(y,z)$ to the generator $g$ as $n\To\infty$ and there exists a random variable $\eta\in {\mathbb L}^1(\F_T)$ such that for each $n\geq 1$ and $t\in\T$,
\begin{equation}
\begin{array}{l}
\Dis\E\left[\left.\sup\limits_{s\in [t,T]}|Y_s^n|^p+\left(\int_t^T|Z_s^n|^2{\rm d}s\right)^{p\over 2}+|K_T^n-K_t^n|^p+\left(\int_t^T|g_n(s,Y_s^n,Z_s^n)|{\rm d}s\right)^p\right|\F_t\right]\\
\ \ \ \ \leq \E\left[\left.\eta\right|\F_t\right],
\end{array}
\end{equation}
then there exists a triple $(Y_\cdot,Z_\cdot,K_\cdot)\in \s^p\times\M^p\times\vcal^{+,p}$ which solves RBSDE (1),
$$\lim\limits_{n\To \infty}\left(\|Y_\cdot^n-Y_\cdot\|_{\s^p}+
\|Z_\cdot^n-Z_\cdot\|_{\M^p}\right)=0,$$
and there exists a subsequence $\{K_\cdot^{n_j}\}$ of $\{K_\cdot^n\}$ such that
$$\lim\limits_{j\To\infty}\sup\limits_{t\in\T}
|K_t^{n_j}-K_t|=0.\vspace{-0.1cm}$$
Furthermore, if $K_\cdot^n$ increases or decrease in $n$, then we have
$$\lim\limits_{n\To \infty} \|K_\cdot^n-K_\cdot\|_{\s^p}=0.$$

{\bf Proof.}\ Since $Y_\cdot^n$ increases or decreases in $n$, there exists a process $Y_\cdot$ such that $\ps$, $Y_t^n\To Y_t$ for each $t\in\T$. In the same way as in the proof of Proposition 3, by Fatou's lemma together with (37) we can deduce that
\begin{equation}
\E\left[\sup\limits_{t\in\T}|Y_t|^p\right]
<+\infty,
\end{equation}
\begin{equation}
\sup\limits_{n\geq 1}|Y^n_t|\leq \left(\E\left[\left.\eta\right|\F_t\right]
\right)^{1\over p},\ \ t\in\T
\end{equation}
and
\begin{equation}
\sup\limits_{n\geq 1}\E\left[\left(\int_0^T|Z_t^n|^2{\rm d}t\right)^{p\over 2}+|K_T^n|^p+\left(\int_0^T|g_n(t,Y_t^n,Z_t^n)|{\rm d}t\right)^p\right]\leq \E[\eta]<+\infty.\vspace{0.1cm}
\end{equation}
For each positive integer $l,q\geq 1$, as in the proof of Proposition 3, we introduce the following two $(\F_t)$-stopping times:
$$
\begin{array}{rll}
\tau_l&:=&\Dis \inf\left\{t\geq 0:\  \left(\E\left[\left.\eta\right|\F_t\right]
\right)^{1\over p}+\int_0^t f_s{\rm d}s\geq l\right\}\wedge T;\\
\sigma_{l,q}&:=&\Dis \inf\left\{t\geq 0:\  \int_0^t \psi_s(l) {\rm d}s\geq q\right\}\wedge \tau_l.
\end{array}
$$
Then we have
\begin{equation}
\mathbb{P}\left(\left\{\omega:\ \exists l_0(\omega), q_0(\omega)\geq 1,\ \RE l\geq l_0(\omega), \RE q\geq q_0(\omega), \ \sigma_{l,q}(\omega)=T\right\}\right)=1.
\end{equation}
Furthermore, since all $g_n$ satisfy (HH) with the same $f_\cdot, \psi_\cdot(r)$ and $\lambda$, and (39) is satisfied, in view of the definitions of $\tau_l$ and $\sigma_{l,q}$, we know that $\as$, for each $l,q,n\geq 1$,
\begin{equation}
\mathbbm{1}_{t\leq \sigma_{l,q}}|g_n(t,Y_t^n,Z_t^n)|\leq \mathbbm{1}_{t\leq \tau_{l}}f_t+\mathbbm{1}_{t\leq \sigma_{l,q}}\psi_t(l)+\lambda |Z_t^n|\vspace{0.1cm}
\end{equation}
with
\begin{equation}
\E\left[\int_0^T \mathbbm{1}_{t\leq \tau_{l}}f_t{\rm d}t\right]\leq l\ \ {\rm and}\ \ \E\left[\int_0^T \mathbbm{1}_{t\leq \sigma_{l,q}}\psi_t(l){\rm d}t\right]\leq q.\vspace{0.2cm}
\end{equation}
The rest proof is divided into 4 steps.\vspace{0.2cm}

{\bf Step 1.}\ We show the convergence of the sequence $\{Y_\cdot^n\}$. For each $n,m\geq 1$, observe that
\begin{equation}
\begin{array}{lll}
\Dis (\bar Y_\cdot,\bar Z_\cdot,\bar V_\cdot)&
:=&\Dis (Y_\cdot^n-Y_\cdot^m,Z_\cdot^n-Z_\cdot^m,\\
&&\Dis \ \ \int_0^\cdot \left(g_n(s,Y_s^n,Z_s^n)-g_m(s,Y_s^m,Z_s^m)\right){\rm d}s+\left(K_\cdot^n-K_\cdot^m\right) )
\end{array}
\end{equation}
satisfies equation (2). It then follows from (ii) of Lemma 1 with $p=2$, $t=0$ and $\tau=\sigma_{l,q}$ that there exists a constant $C>0$ such that for each $n,m,l,q\geq 1$,
\begin{equation}
\begin{array}{ll}
&\Dis\E\left[\sup\limits_{t\in [0,T]} |Y_{t\wedge\sigma_{l,q}}^n
-Y_{t\wedge\sigma_{l,q}}^m|^2\right]\\
\leq &\Dis C\E\left[|Y_{\sigma_{l,q}}^n-Y_{\sigma_{l,q}}^m|^2 +\sup\limits_{t\in [0,T]}\left(\int_{t\wedge \sigma_{l,q}}^{\sigma_{l,q}}(Y^n_s-Y_s^m)
\left({\rm d}K_s^n-{\rm d}K_s^m\right) \right)^+\right.\vspace{0.1cm}\\
&\hspace{1cm}\Dis +\left. \int_{0}^{\sigma_{l,q}}|Y^n_t-Y_t^m|
\left|g_n(t,Y_t^n,Z_t^n)-g_m(t,Y_t^m,Z_t^m)\right| {\rm d}t\right].
\end{array}
\end{equation}
Furthermore, note that $Y_t^n\geq L_t$ for each $t\in\T$ and $n\geq 1$ and that $\int_0^T(Y_t^n-L_t){\rm d}K_t^n=0$ for each $n\geq 1$. It follows that for each $t\in \T$ and $l,q,m,n\geq 1$,
\begin{equation}
\begin{array}{ll}
&\Dis \int_{t\wedge \sigma_{l,q}}^{\sigma_{l,q}}(Y^n_s-Y_s^m)
\left({\rm d}K_s^n-{\rm d}K_s^m\right)\\
= &\Dis \int_{t\wedge \sigma_{l,q}}^{\sigma_{l,q}}\left[(Y^n_s-L_s)
-(Y_s^m-L_s)\right]
{\rm d}K_s^n-\int_{t\wedge \sigma_{l,q}}^{\sigma_{l,q}}\left[(Y^n_s-L_s)
-(Y_s^m-L_s)\right]
{\rm d}K_s^m\\
\leq &\Dis \int_{t\wedge \sigma_{l,q}}^{\sigma_{l,q}}(Y_s^n-L_s)
{\rm d}K_s^n+\int_{t\wedge \sigma_{l,q}}^{\sigma_{l,q}}(Y_s^m-L_s){\rm d}K_s^m=0.
\end{array}
\end{equation}
Combining (42), (45) and (46) together with H\"{o}lder's inequality yields that
\begin{equation}
\begin{array}{ll}
&\Dis\E\left[\sup\limits_{t\in [0,T]} |Y_{t\wedge\sigma_{l,q}}^n
-Y_{t\wedge\sigma_{l,q}}^m|^2\right]\\
\leq &\Dis C\E\left[|Y_{\sigma_{l,q}}^n-Y_{\sigma_{l,q}}^m|^2 +2\int_0^T |Y^n_t-Y_t^m|\left(\mathbbm{1}_{t\leq \tau_l}f_t+\mathbbm{1}_{t\leq\sigma_{l,q}}
\psi_t(l)\right){\rm d}t\right]\vspace{0.1cm}\\
&\Dis+2C\lambda\left(\E\left[\left(\int_{0}^{
\sigma_{l,q}}|Y^n_t-Y_t^m|^2{\rm d}t\right)^{p\over 2(p-1)}\right]\right)^{p-1\over p}\\
&\Dis \hspace{1cm}\times \left(\E\left[\left(\int_{0}^T\left(|Z_t^n|+
|Z_t^m|\right)^2{\rm d}t\right)^{p\over 2}\right]\right)^{1\over p}.
\end{array}
\end{equation}
Thus, note that $Y_t^n\To Y_t$ for each $t\in\T$. In view of the definitions of $\tau_l$ and $\sigma_{l,q}$, (39), (40) and (43), by (47) and Lebesgue's dominated convergence theorem we can deduce that for each $l,q\geq 1$, as $n,m\To\infty$,
$$\E\left[\sup\limits_{t\in [0,T]} |Y_{t\wedge\sigma_{l,q}}^n
-Y_{t\wedge\sigma_{l,q}}^m|^2\right]\To 0,$$
which implies that for each $l,q\geq 1$, as $n,m\To\infty$,
$$\sup\limits_{t\in [0,T]} |Y_{t\wedge\sigma_{l,q}}^n
-Y_{t\wedge\sigma_{l,q}}^m|\To 0\ {\rm in\ probability}\ \mathbb{P}.$$
And, by (41) and the monotonicity of $Y_\cdot^n$ with respect to $n$ we know that $\ps$,
\begin{equation}
\sup\limits_{t\in [0,T]} |Y_t^n
-Y_t|\To 0,\ \ {\rm as}\ n\To\infty.
\end{equation}
So, $Y_\cdot$ is a continuous process. Finally, note that $|Y_\cdot^n|\leq |Y_\cdot^1|+|Y_\cdot|$ for each $n\geq 1$ and that (38) is satisfied. From (48) and Lebesgue's dominated convergence theorem it follows that
\begin{equation}
\lim\limits_{n\To\infty}\|Y_\cdot^n-Y_\cdot\|_{\s^p}^p =\lim\limits_{n\To\infty}\E\left[\sup\limits_{t\in [0,T]} |Y_t^n-Y_t|^p\right]=0.
\end{equation}

{\bf Step 2.}\ We show the convergence of the sequence $\{Z_\cdot^n\}$. Note that (44) solves (2). It follows from (i) of Lemma 1 with $t=0$ and $\tau=T$ that there exists a constant $C'>0$ such that for each $m,n\geq 1$,
$$
\begin{array}{ll}
&\Dis\E\left[\left(\int_0^T|Z_t^n-Z_t^m|^2{\rm d}t\right)^{p\over 2}\right]\\
\leq & \Dis C'\E\left[\sup\limits_{t\in [0,T]}|Y_t^n-Y_t^m|^p+\sup\limits_{t\in [0,T]}\left[\left(\int_t^T (Y_s^n-Y_s^m)\left({\rm d}K_s^n-{\rm d}K_s^m\right)\right)^+\right]^{p\over 2}\right]\\
&\Dis + C'\E\left[\left(\int_{0}^T|Y^n_t-Y_t^m|
\left|g_n(t,Y_t^n,Z_t^n)-g_m(t,Y_t^m,Z_t^m)\right| {\rm d}t\right)^{p\over 2}\right].
\end{array}
$$
Then, in view of (46), it follows from H\"{o}lder's inequality that for each $m,n\geq 1$,
$$
\begin{array}{ll}
&\Dis\E\left[\left(\int_0^T|Z_t^n-Z_t^m|^2{\rm d}t\right)^{p\over 2}\right]\\
\leq &\Dis C'\E\left[\sup\limits_{t\in [0,T]}|Y_t^n-Y_t^m|^p\right]+C'\left(\E\left[\sup\limits_{t\in [0,T]}|Y_t^n-Y_t^m|^p\right]\right)^{1\over 2}\\
&\Dis\ \ \ \ \  \cdot\left(\E\left[\left(\int_{0}^T
\left(|g_n(t,Y_t^n,Z_t^n)|+|g_n(t,Y_t^m,Z_t^m)|\right) {\rm d}t\right)^{p}\right]\right)^{1\over 2},
\end{array}
$$
from which together with (49) and (40) yields that there exists a process $Z_\cdot\in\M^p$ such that
\begin{equation}
\lim\limits_{n\To\infty}\|Z_\cdot^n-Z_\cdot\|_{\M^p}^p
=\lim\limits_{n\To\infty}\E\left[\left(\int_0^T|Z_t^n
-Z_t|^2{\rm d}t\right)^{p\over 2}\right]=0.
\end{equation}

{\bf Step 3.}\ We show the convergence of the sequence $\{K_\cdot^n\}$. Since $g_n$ tends locally uniformly in $(y,z)$ to the generator $g$ as $n\To+\infty$, by (48), (50), (40), (42) and (43) we can deduce that there exists a subsequence $\{n_j\}$ of $\{n\}$ such that for each $l,q\geq 1$,
$$\lim\limits_{j\To \infty}\int_0^{\sigma_{l,q}}|g_{n_j}(t,Y_t^{n_j},Z_t^{n_j})-
g(t,Y_t,Z_t)|{\rm d}t=0.$$
Then, in view of (41), we have
\begin{equation}
\lim\limits_{j\To \infty}\sup\limits_{t\in\T}\left|\int_0^t
g_{n_j}(t,Y_t^{n_j},Z_t^{n_j}){\rm d}t-\int_0^t
g(t,Y_t,Z_t){\rm d}t\right|=0.
\end{equation}
Combining (48), (50) and (51) yields that $\ps$, for each $t\in\T$,
$$
K_t^{n_j}=Y_0^{n_j}-Y_t^{n_j}-\int_0^t
g_{n_j}(s,Y_s^{n_j},Z_s^{n_j}){\rm d}s-\int_0^t{\rm d}V_s+\int_0^tZ_s^{n_j}{\rm d}B_s
$$
tends to
$$
K_t:=Y_0-Y_t-\int_0^t
g(s,Y_s,Z_s){\rm d}s-\int_0^t{\rm d}V_s+\int_0^tZ_s{\rm d}B_s
$$
as $j\To \infty$ and that
\begin{equation}
\lim\limits_{j\To\infty}\sup\limits_{t\in\T}
|K_t^{n_j}-K_t|=0.
\end{equation}
Hence, $K_\cdot$ is a continuous process.\vspace{0.2cm}

{\bf Step 4.}\ We show that $K_\cdot\in\vcal^{+,p}$ and $(Y_\cdot,Z_\cdot,
K_\cdot)\in \s^p\times\M^p\times\vcal^{+,p}$ is a solution of RBSDE (1). In fact, by Fatou's lemma with (52) and (40) we get that
$$
\begin{array}{lll}
\Dis\E\left[\sup\limits_{t\in\T}|K_t|^p\right]
&=&\E\left[\lim\limits_{j\To\infty}
\sup\limits_{t\in\T}|K_t^{n_j}|^p\right]\leq \Dis \liminf\limits_{j\To\infty}
\E\left[\sup\limits_{t\in\T}|K_t^{n_j}|^p\right]
\vspace{0.1cm}\\
&\leq& \Dis \sup\limits_{j\geq 1}\E\left[|K_T^{n_j}|^p\right]\leq  \E\left[\eta\right]<+\infty.
\end{array}
$$
Hence, $K_\cdot\in\vcal^{+,p}$ and $(Y_\cdot,Z_\cdot,K_\cdot)\in \s^p\times\M^p\times\vcal^{+,p}$ solves
$$
Y_t=\xi+\int_t^Tg(s,Y_s,Z_s){\rm d}s+\int_t^T{\rm d}V_s+\int_t^T{\rm d}K_s-\int_t^TZ_s {\rm d}B_s,\ \ t\in\T.
$$
Since $Y_t^n\geq L_t,\ n\geq 1$ and $Y_t^n\To Y_t$ for each $t\in\T$, we have $Y_t\geq L_t$ for each $t\in\T$. Furthermore, in view of (48) and (52), it follows that
$$
\int_0^T(Y_t-L_t){\rm d}K_t=\lim\limits_{j\To \infty}\int_0^T(Y_t^{n_j}-L_t){\rm d}K_t^{n_j}=0.
$$
So, $(Y_\cdot,Z_\cdot,
K_\cdot)\in \s^p\times\M^p\times\vcal^{+,p}$ solves RBSDE (1).

Finally, if $K_\cdot^n$ increases or decrease in $n$, then $\ps$, for each $t\in\T$, $K_t^n\To K_t$ as $n\To\infty$ and $|K_t^n|\leq |K^1_t|+|K_t|$. Thus, it follows from Dini's theorem and Lebesgue's dominated convergence theorem that
$$\lim\limits_{n\To \infty} \|K_\cdot^n-K_\cdot\|_{\s^p}=0.$$
Proposition 4 is then proved.\vspace{0.2cm}\hfill $\Box$

{\bf Remark 4}\ In the case when $L_\cdot=-\infty$ and $K^n_\cdot=0$ for each $n\geq 1$, by Proposition 4 we can get the approximation result for $L^p$ solutions of non-reflected BSDEs.\vspace{0.2cm}

{\bf Proposition 5} (Comparison theorem) Let $p>1$, $V_\cdot^i\in \vcal^p$, $g^i$ be a generator, $\xi^i$ and $L^i_\cdot$ satisfy (H5), and $(Y_\cdot^i,Z_\cdot^i,K_\cdot^i)\in \s^p\times \M^p\times \vcal^{+,p}$ be a solution of RBSDE $(\xi^i,g^i+{\rm d}V^i,L^i)$ for each $i=1,2$. If $\xi^1\leq \xi^2$, ${\rm d}V^1_\cdot\leq {\rm d}V^2_\cdot$, $L^1_\cdot\leq L^2_\cdot$, and either
$$
\left\{
\begin{array}{l}
g^1\ {\rm satisfies\ (H1)\ and \ (H2)};\\
\as,\ \ \mathbbm{1}_{\{Y^1_t>Y^2_t\}}
\left(g^1(t,Y_t^2,Z_t^2)
-g^2(t,Y_t^2,Z_t^2)\right)\leq 0
\end{array}
\right.
$$
or
$$
\left\{
\begin{array}{l}
g^2\ {\rm satisfies\ (H1)\ and \ (H2)};\\
\as,\ \ \mathbbm{1}_{\{Y^1_t>Y^2_t\}}
\left(g^1(t,Y_t^1,Z_t^1)
-g^2(t,Y_t^1,Z_t^1)\right)\leq 0
\end{array}
\right.
$$
is satisfied, then $\ps$, $Y_t^1\leq Y_t^2$ for each $t\in \T$.\vspace{0.2cm}

{\bf Proof.}\ By It\^{o}-Tanaka's formula we know that for each $t\in\T$,
$$
\begin{array}{lll}
\Dis  (Y^1_t-Y^2_t)^+ &\leq &\Dis (\xi^1-\xi^2)^++\int_t^T {\rm sgn}((Y^1_s-Y^2_s)^+)({\rm d}V_s^1-{\rm d}V_s^2)\\
&&\Dis +\int_t^T{\rm sgn}((Y^1_s-Y^2_s)^+)\left(g^1(s,Y_s^1,Z_s^1)
-g^2(s,Y_s^2,Z_s^2)\right){\rm d}s\\
&&\Dis+\int_t^T {\rm sgn}((Y^1_s-Y^2_s)^+)\left({\rm d}K_s^1-{\rm d}K_s^2\right)\\
&&\Dis +\int_t^T {\rm sgn}((Y^1_s-Y^2_s)^+)(Z_s^1-Z_s^2){\rm d}B_s.
\end{array}
$$
Since $L_t^1\leq L_t^2\leq Y_t^2$, $L_t^1\leq Y_t^1$, $t\in \T$ and $\int_0^T(Y^1_s-L_s^1){\rm d}K_s^1=0$, we have
$$
\begin{array}{lll}
\Dis \int_t^T {\rm sgn}((Y^1_s-Y^2_s)^+)\left({\rm d}K_s^1-{\rm d}K_s^2\right)\\
\ \ \leq  \Dis \int_t^T {\rm sgn}((Y^1_s-Y^2_s)^+){\rm d}K_s^1\leq \int_t^T {\rm sgn}((Y^1_s-L_s^1)^+){\rm d}K_s^1\\
\ \ \Dis =\int_t^T \mathbbm{1}_{\{Y^1_s>L_s^1\}}|Y^1_s-L_s^1|^{-1}(Y^1_s-L_s^1){\rm d}K_s^1=0.
\end{array}
$$
Thus, noticing that $\xi^1\leq \xi^2$ and ${\rm d}V^1_t\leq {\rm d}V^2_t$ for each $t\in\T$, we can get that
$$
\begin{array}{lll}
\Dis  (Y^1_t-Y^2_t)^+ &\leq &\Dis \int_t^T{\rm sgn}((Y^1_s-Y^2_s)^+)\left(g^1(s,Y_s^1,Z_s^1)
-g^2(s,Y_s^2,Z_s^2)\right){\rm d}s\\
&&\Dis +\int_t^T {\rm sgn}((Y^1_s-Y^2_s)^+)(Z_s^1-Z_s^2){\rm d}B_s,\ \ t\in\T.
\end{array}
$$
Now, in view of the assumptions of $g^1$ and $g^2$, the rest proof runs as the proof of Theorem 1 in \citet{Fan12}. The only difference lies in that in order to deal with the $L^p$ solution we need to use
$$\E\left[|XY|\right]\leq \left(\E\left[|X|^p\right]\right)^{1\over p} \left(\E\left[|Y|^{p\over p-1}\right]\right)^{p-1\over p}$$
instead of the inequality
$$\E\left[|XY|\right]\leq \left(\E\left[|X|^2\right]\right)^{1\over 2} \left(\E\left[|Y|^2\right]\right)^{1\over 2}\vspace{0.1cm}$$
for any $\F_T$-measurable random variables $X$ and $Y$.\vspace{0.2cm} \hfill $\Box$

From Proposition 5, the following corollary is immediate. \vspace{0.2cm}

{\bf Corollary 1}\ Let $p>1$, $V_\cdot^i\in \vcal^p$, $g^i$ be a generator, $\xi^i$ and $L^i_\cdot$ satisfy (H5), and $(Y_\cdot^i,Z_\cdot^i,K_\cdot^i)\in \s^p\times \M^p\times \vcal^{+,p}$ be a solution of RBSDE $(\xi^i,g^i+{\rm d}V^i,L^i)$ for each $i=1,2$. If $\xi^1\leq \xi^2$, ${\rm d}V^1_\cdot\leq {\rm d}V^2_\cdot$, $L^1_\cdot\leq L^2_\cdot$, $g^1$ or $g^2$ satisfies (H1) and (H2), and
$$\as,\ \ g^1(t,y,z)\leq g^2(t,y,z)$$
for each $(y,z)\in \R\times \R^d$, then $\ps$, $Y_t^1\leq Y_t^2$ for each $t\in \T$.\vspace{0.2cm}

{\bf Remark 5}\ In the case when $L^1_\cdot=L^2_\cdot=-\infty$ and $K^1_\cdot=K^2_\cdot=0$, by Proposition 5 and Corollary 1 we can get the comparison result for $L^p$ solutions of non-reflected BSDEs. In addition, it should be noted that Proposition 5 and Corollary 1 improves greatly the corresponding results obtained in \citet{El97b,Fan12,Ham12,Kli12,Lep05,Roz12} and etc.

\section{Existence, uniqueness and approximation results}

In this section, based on the results obtained in previous sections, we will establish some existence, uniqueness and approximation results on $L^p$ solutions of BSDEs and RBSDEs with $L^p\ (p>1)$ data under weaker assumptions, which answers those questions put forward in the Introduction.\vspace{-0.1cm}

\subsection{Non-reflected BSDEs}\vspace{0.1cm}

Let us start with the following existence and uniqueness result---Proposition 6. It improves Corollary 2 of \citet{Fan15} in the one-dimensional case, where $V_\cdot\equiv 0$ and the $\varphi_\cdot(r)$ in (H3) is assumed to be in $\hcal^1$.\vspace{0.2cm}

{\bf Proposition 6}\ Let $p>1$, $V_\cdot\in\vcal^{p}$ and let $g$ satisfy assumptions (H1), (H2s), (H3) and (H4w). Then for each $\xi\in\Lp$, the following equation, denoted by BSDE $(\xi,g+{\rm d}V)$ here and hereafter,
$$
Y_t=\xi+\int_t^T g(s,Y_s,Z_s){\rm d}s+\int_t^T {\rm d}V_s-\int_t^T Z_s{\rm d}B_s,\ \ t\in \T
$$
admits a unique solution in $\s^p\times\M^p$.\vspace{0.2cm}

{\bf Proof.} Note that (H2s)$\Longrightarrow$(H2). The uniqueness part follows immediately from Proposition 5 or Corollary 1 in view of Remark 5. In the sequel, we prove the existence part. Let $p>1$, $\xi\in\Lp$, $V_\cdot\in\vcal^{p}$ and let the generator $g$ satisfy (H1) with $\rho(\cdot)$, (H2s) with $\lambda$, (H3) with $\varphi_\cdot(r)$ and (H4w).

We first assume that $g$ is bounded.
It then follows from Corollary 2 in \citet{Fan15} that the following BSDE
$$
\bar Y_s=\xi+V_T+\int_t^T g(s,\bar Y_s-V_s,\bar Z_t){\rm d}s-\int_t^T\bar Z_s{\rm d}B_s,\ \ t\in\T
$$
admits a unique solution $(\bar Y_\cdot,\bar Z_\cdot)\in \s^p\times\M^p$. Then the pair $(Y_\cdot,Z_\cdot):=(\bar Y_\cdot-V_\cdot,\bar Z_\cdot)$ is the unique solution of BSDE $(\xi,g+{\rm d}V)$ in $\s^p\times\M^p$.

Now suppose that $g$ is bounded from below. Write $g_n=g\wedge n$. Then $g_n$ is bounded, nondecreasing in $n$ and tends locally uniformly to $g$ as $n\To\infty$, and it is not difficult to check that all $g_n$ satisfy (H1), (H2s) and (H3) with the same $\rho(\cdot)$, $\lambda$ and $\varphi_\cdot(r)$ as well as (H4w). Then by the first step of the proof there exists a unique solution $(Y_\cdot^n,Z_\cdot^n)$ of BSDE $(\xi,g_n+{\rm d}V)$. By Corollary 1 together with Remark 5, $Y_\cdot^n$ increases in $n$. Furthermore, $\as$, for each $n\geq 1$ and $(y,z)\in\R\times\R^d$,\vspace{-0.2cm}
$$
\begin{array}{lll}
\Dis |g_n(\cdot,y,z)|&\leq &\Dis |g_n(\cdot,y,z)-g_n(\cdot,y,0)|+|g_n(\cdot,y,0)-g_n(\cdot,0,0)|
+|g_n(\cdot,0,0)|\\
&\leq & \Dis \lambda |z|+\varphi_\cdot(|y|)+|g(\cdot,0,0)|,
\end{array}
$$
which means that all $g_n$ satisfy (HH) with the same
$$
f_\cdot:=|g(\cdot,0,0)|,\ \  \psi_\cdot(r):=\varphi_\cdot(r)\ \  {\rm and}\ \  \lambda.
$$
It then follows from Proposition 1 that (37) in Proposition 4 holds with
$$\eta:=C\left[|\xi|^p+|V|_T^p+\left(\int_0^T |g(t,0,0)|{\rm d}t\right)^p+1 \right],$$
where $K_\cdot^n\equiv 0$ and $C>0$ is a constant depending only $p,\lambda,A,T$. Thus, we have checked all the conditions in Proposition 4, and then in view of Remark 4, it follows from Proposition 4 that BSDE $(\xi,g+{\rm d}V)$ admits a solution in $\s^p\times\M^p$.

Finally, in the general case, we can approximate $g$ by the sequence $g_n$, where $g_n:=g\vee (-n),\ n\geq 1$. By the previous step there exists a unique solution $(Y^n_\cdot,Z^n_\cdot)\in \s^p\times\M^p$ of BSDE $(\xi,g_n+{\rm d}V)$ for each $n\geq 1$. Repeating arguments in the proof of the previous step yields that $(Y^n_\cdot,Z^n_\cdot)$ converges in $\s^p\times \M^p$ to the unique solution of BSDE $(\xi,g+{d}V)$.\vspace{0.2cm}\hfill $\Box$

By Propositions 1 and 4-6, we can prove the following Theorems 1 and 2, which further extend Proposition 6 to the case of BSDEs with generators satisfying the weaker assumptions (H2) or (H2w) than (H2s). They generalize respectively Theorems 2-3 in \citet{Fan12} and Theorem 4.1 in \citet{Bri07}, where some stronger assumptions are assumed to be satisfied.\vspace{0.2cm}

{\bf Theorem 1}\ Let $p>1$, $V_\cdot\in\vcal^{p}$ and let $g$ satisfy (H1), (H2), (H3) and (H4w). Then for each $\xi\in\Lp$, BSDE $(\xi,g+{\rm d}V)$ admits a unique solution in $\s^p\times\M^p$.\vspace{0.2cm}

{\bf Proof.} The uniqueness part is a direct consequence of Proposition 5 in view of Remark 5.

Now, we prove the existence part. Firstly, let $p>1$, $\xi\in\Lp$, $V_\cdot\in\vcal^{p}$ and let the generator $g$ satisfy assumptions (H1) with $\rho(\cdot)$, (H2) with $\phi(\cdot)$, (H3) with $\varphi_\cdot(r)$ and (H4w). By a similar argument to that in the proof of the existence part of Theorem 1 of \citet{Ma13} we can prove that for each $n\geq 1$ and $(y,z)\in \R\times\R^d$, the following function
$$
g_n(\omega,t,y,z):=\inf\limits_{u\in\R^d}
\left[g(\omega,t,y,u)+(n+2A)|u-z|\right]
$$
is well defined and $(\F_t)$-progressively measurable, $\as$, $g_n$ increases in $n$ and  converges uniformly in $(y,z)$ to the generator $g$,
$$
\sup\limits_{n\geq 1}|g_n(\cdot,0,0)|\leq |g(\cdot,0,0)|+\phi(2A)
$$
and all $g_n$ satisfy (H1) with the same $\rho(\cdot)$, (H2) with the same $\phi(\cdot)$, (H3) with the same $\varphi_\cdot(r)+2\phi(2A)$, (H4) and (H2s) with $\lambda:=n+2A$. It then follows from Proposition 6 that there exists a unique solution $(Y_\cdot^n,Z_\cdot^n)\in
\s^p\times\M^p$ of BSDE $(\xi,g_n+{\rm d}V)$ for each $n\geq 1$. By Corollary 1 together with Remark 5, $Y_\cdot^n$ increases in $n$. Furthermore, it follows from Remark 1 that $\as$, for each $n\geq 1$ and $(y,z)\in\R\times\R^d$,\vspace{-0.2cm}
$$
\begin{array}{lll}
\Dis |g_n(\cdot,y,z)|&\leq &\Dis |g_n(\cdot,y,z)-g_n(\cdot,y,0)|+|g_n(\cdot,y,0)-g_n(\cdot,0,0)|
+|g_n(\cdot,0,0)|\\
&\leq & \Dis A|z|+A+\varphi_\cdot(|y|)+2\phi(2A)+|g(\cdot,0,0)|+\phi(2A),
\end{array}
$$
which means that all $g_n$ satisfy (HH) with the same
$$
f_\cdot:=|g(\cdot,0,0)|+3\phi(2A)+A,\ \  \psi_\cdot(r):=\varphi_\cdot(r)\ \  {\rm and}\ \  \lambda:=A.
$$
It then follows Proposition 1 that (37) in Proposition 4 holds with
$$\eta:=C\left[|\xi|^p+|V|_T^p+\left(\int_0^T |g(t,0,0)|{\rm d}t\right)^p+1 \right],$$
where $K_\cdot^n\equiv 0$ and $C>0$ is a constant depending only $p,A,T$. Thus, we have checked all the conditions in Proposition 4, and then in view of Remark 4, the conclusion of Theorem 1 follows from Proposition 4 immediately.
\vspace{0.2cm}\hfill$\Box$

{\bf Theorem 2}\ Let $p>1$, $V_\cdot\in\vcal^{p}$ and let $g$ satisfy (H1), (HH) and (H4s). Then for each $\xi\in\Lp$, BSDE $(\xi,g+{\rm d}V)$ admits a maximal (resp. minimal) solution $(Y_\cdot,Z_\cdot)$ in $\s^p\times\M^p$, i.e., if $(Y'_\cdot,Z'_\cdot)$ is also a solution of BSDE $(\xi,g+{\rm d}V)$ in $\s^p\times\M^p$, then $\ps$, $Y_t\geq Y'_t$ (resp. $Y_t\leq Y'_t$) for each $t\in \T$.\vspace{0.2cm}

{\bf Proof.} We only prove the case of the maximal solution. In the same way, we can prove another case.

Assume that $p>1$, $\xi\in\Lp$, $V_\cdot\in\vcal^{p}$ and let the generator $g$ satisfy assumptions (H1) with $\rho(\cdot)$, (HH) with $f_\cdot$, $\varphi_\cdot(r)$ and $\lambda$ as well as (H4s). It is not very hard to prove that for each $n\geq 1$ and $(y,z)\in \R\times\R^d$, the following function
\begin{equation}
g_n(\omega,t,y,z):=\sup\limits_{u\in\R^d}
\left[g(\omega,t,y,u)-(n+2\lambda)|u-z|\right]
\end{equation}
is well defined and $(\F_t)$-progressively measurable, $\as$, $g_n$ decreases in $n$ and  converges locally uniformly in $(y,z)$ to the generator $g$ as $n\To \infty$, and all $g_n$ satisfy (H1) with the same $\rho(\cdot)$, (HH) with the same $f_\cdot$, $\varphi_\cdot(r)$ and $\lambda$, (H4) and (H2s) with $\lambda:=n+2\lambda$. Note by (i) and (ii) of Remark 2 that (HH)$\Rightarrow$(H3) and (H4)$\Rightarrow$(H4w). It then follows from Proposition 6 that there exists a unique solution $(Y_\cdot^n,Z_\cdot^n)\in \s^p\times\M^p$ of BSDE $(\xi,g_n+{\rm d}V)$ for each $n\geq 1$. By Corollary 1 together with Remark 5, $Y_\cdot^n$ decreases in $n$. Furthermore, it follows from Proposition 1 that (37) in Proposition 4 holds true with
$$\eta:=C\left[|\xi|^p+|V|_T^p+\left(\int_0^T f_t{\rm d}t\right)^p+1\right],$$
where $K_\cdot^n\equiv 0$ and $C>0$ is a constant depending only $p,\lambda,A,T$. Thus, we have checked all the conditions in Proposition 4, and then in view of Remark 4, it follows from Proposition 4 that BSDE $(\xi,g+{\rm d}V)$ admits a solution $(Y_\cdot,Z_\cdot)$ in $\s^p\times\M^p$ such that
\begin{equation}
\lim\limits_{n\To\infty}(\|Y_\cdot^n
-Y_\cdot\|_{\s^p}+\|Z_\cdot^n
-Z_\cdot\|_{\M^p})=0.
\end{equation}

Finally, we show that $(Y_\cdot,Z_\cdot)$ is just the maximal solution of BSDE $(\xi,g+{\rm d}V)$ in $\s^p\times\M^p$. In fact, if $(Y'_\cdot,Z'_\cdot)$ is also a solution of BSDE $(\xi,g+{\rm d}V)$ in $\s^p\times\M^p$, then noticing that for each $n\geq 1$, $g_n\geq g$ and $g_n$ satisfies (H1) and (H2) due to (H2s)$\Rightarrow$(H2), it follows from Corollary 1 together with Remark 5 that $Y_t^n\geq Y'_t$ for each $t\in\T$ and $n\geq 1$. Thus, by (54) we know that $\ps$, $Y_t\geq Y'_t$ for each $t\in \T$. Theorem 2 is then proved.\vspace{0.2cm}\hfill$\Box$

In view of (ii) of Remark 2, the following corollary follows from Theorem 2 immediately. \vspace{0.2cm}

{\bf Corollary 2}\ Let $p>1$, $V_\cdot\in\vcal^{p}$ and let $g$ satisfy (H1), (H2w), (H3) and (H4s). Then for each $\xi\in\Lp$, BSDE $(\xi,g+{\rm d}V)$ admits a maximal (resp. minimal) solution $(Y_\cdot,Z_\cdot)$ in $\s^p\times\M^p$.\vspace{0.2cm}

By Corollary 1 together with Remark 5 and the proof of Theorem 2 it is easy to verify that under assumptions (H1), (HH) and (H4s), the comparison theorem for the maximal (resp. minimal) solutions of BSDEs holds. More precisely, we have

{\bf Corollary 3}\ Let $p>1$ and for $i=1,2$, assume that $\xi^i\in\Lp$, $V_\cdot^i\in \vcal^p$, $g^i$ satisfies (H1), (HH) and (H4s), and that $(Y_\cdot^i,Z_\cdot^i)\in \s^p\times \M^p$ is the maximal (resp. minimal) solution of BSDE $(\xi^i,g^i+{\rm d}V^i)$. If $\xi^1\leq \xi^2$, ${\rm d}V^1_\cdot\leq {\rm d}V^2_\cdot$, and
$$\as,\ \ g^1(t,y,z)\leq g^2(t,y,z)$$
for each $(y,z)\in \R\times \R^d$, then $\ps$, $Y_t^1\leq Y_t^2$ for each $t\in \T$.\vspace{0.1cm}

\subsection{Reflected BSDEs}\vspace{0.2cm}

The following theorem shows that under conditions of (H1), (H2w) and (H5) with $g(\cdot,0,0)\in \hcal^p$, (H6) is necessary to ensure the existence of $L^p$ solutions for RBSDEs with $L^p\ (p>1)$ data, which is one of our main results.\vspace{0.2cm}

{\bf Theorem 3}\ Assume that $p>1$, $V_\cdot\in \vcal^p$, the generator $g$ satisfies (H1) and (H2w) with $g(\cdot,0,0)\in \hcal^p$, and that (H5) holds for $\xi$ and $L_\cdot$. If RBSDE (1) admits a solution $(Y_\cdot,Z_\cdot,K_\cdot)\in \s^p\times \M^p\times \vcal^{+,p}$, then $g(\cdot,Y_\cdot,0)\in \hcal^p$. So (H6) holds.\vspace{0.2cm}

{\bf Proof.}\ By (iv) of Remark 2 we know that $g$ satisfies (A) with $\bar f_\cdot:=|g(\cdot,0,0)|+f_\cdot+A$, $\bar\mu:=\mu+A$
and $\bar\lambda:=\lambda$. It then follows from (4) in Lemma 3 that
$$
\E\left[\left(\int_0^T |g(t,Y_t,Z_t)|{\rm d}t\right)^p\right]<+\infty.
$$
Then, by (H2w) together with H\"{o}lder's inequality we can deduce that
$$
\begin{array}{ll}
&\Dis\E\left[\left(\int_0^T |g(t,Y_t,0)|{\rm d}t\right)^p\right]\vspace{0.1cm}\\
\leq &\Dis
4^p\E\left[\left(\int_0^T |g(t,Y_t,Z_t)|{\rm d}t\right)^p\right]+4^p\E\left[\left(\int_0^T f_t{\rm d}t\right)^p\right]\vspace{0.1cm}\\
&\Dis +(4\mu T)^p\E\left[\sup\limits_{t\in\T}|Y_t|^p\right]+(4\lambda)^p \E\left[\left(\int_0^T |Z_t|^2{\rm d}t\right)^{p\over 2}\right]\\
<&+\infty.
\end{array}
$$
Thus, Theorem 3 is proved. \vspace{0.2cm} $\Box$ \hfill

The following Proposition 7 establishes an a priori estimate on the $L^p$ solution of the penalization equation of RBSDE $(\xi,g+{\rm d}V,L_\cdot)$ with $L^p\ (p>1)$ data under the assumptions of (H1), (H2) or (H2w), (H3), (H5) and (H6), which will play an important role in the proofs of the following Theorems 4 and 5. The proof is based on Propositions 1-2 and 5, Theorems 1-2, and Corollaries 1-3.\vspace{0.2cm}

{\bf Proposition 7}\ Let $p>1$, $V_\cdot\in\vcal^{p}$ and let $g$ satisfy (H1), (H2w) (resp. (H2)), (H3) and (H4s) (resp. (H4w)). Assume that (H5) and (H6) hold for $\xi$, $L_\cdot$ and some $X_\cdot$. For each $n\geq 1$, let $(Y_\cdot^n,Z_\cdot^n)\in \s^p\times\M^p$ be the maximal or minimal (resp. unique) solution of the penalization equation (16) with (17) (Recall Corollary 2 (resp. Theorem 1)). Then, (18) appearing in Proposition 3 holds true.\vspace{0.2cm}

{\bf Proof.}\ We only prove the case that (H2w) and (H4s) are satisfied and $(Y_\cdot^n,Z_\cdot^n)$ is the maximal solution. In the same way, we can prove other cases.

We first show that (C) appearing in Proposition 2 holds true for $Y^n_\cdot$ and $g$. In fact, since $X_\cdot\in\M^p+\vcal^p$ and the Brownian filtration has the representation property, there exist $H_\cdot\in\M^p$ and $C_\cdot\in\vcal^p$ such that
\begin{equation}
X_t=X_T-\int_t^T{\rm d}C_s-\int_t^TH_s{\rm d}B_s,\ \ t\in\T.
\end{equation}
It follows from (H2w) that $\as$,
$$
|g(\cdot,X_\cdot,H_\cdot)|\leq |g(\cdot,X_\cdot,0)|+f_\cdot+\mu |X_\cdot|+\lambda |H_\cdot|,
$$
from which together with (H6) we know that $g(\cdot,X_\cdot,H_\cdot)\in \hcal^p$.
Then, the equation (55) can be rewritten in the form
$$
\begin{array}{lll}
X_t&=&\Dis X_T+\int_t^Tg(s,X_s,H_s){\rm ds}+\int_t^T{\rm d}V_s-\int_t^T\left(g^+(s,X_s,H_s){\rm d}s+{\rm d}C^+_s+{\rm d}V^+_s\right)\\
&&\Dis +\int_t^T\left(g^-(s,X_s,H_s){\rm d}s+{\rm d}C^-_s+{\rm d}V^-_s\right)-\int_t^TH_s{\rm d}B_s,\ \ t\in\T.
\end{array}
$$
On the other hand, it follows from Corollary 2 that there exists a maximal solution $(\bar X_\cdot,\bar Z_\cdot)\in \s^p\times\M^p$ of the BSDE
$$
\begin{array}{lll}
\bar X_t&=&\Dis X_T\vee \xi+\int_t^Tg(s,\bar X_s,\bar Z_s){\rm ds}+\int_t^T{\rm d}V_s\\
&&\Dis +\int_t^T\left(g^-(s,X_s,H_s){\rm d}s+{\rm d}C^-_s+{\rm d}V^-_s\right)-\int_t^T\bar Z_s{\rm d}B_s,\ \ t\in\T.
\end{array}
$$
Note that $g$ satisfies (HH) by (ii) of Remark 2 . It then follows from Proposition 1 with $g_n\equiv g$ that
$$
g(\cdot,\bar X_\cdot,\bar Z_\cdot)\in \hcal^p,\vspace{0.1cm}
$$
and by (H2w) and Remark 1 together with H\"{o}lder's inequality,
$$g(\cdot,\bar X_\cdot,0)\in \hcal^p.$$
Furthermore, it follows from Corollary 3 together with (H6) that $L_t\leq X_t\leq \bar X_t$ for each $t\in \T$. Therefore, we have
$$
\begin{array}{lll}
\bar X_t&=&\Dis X_T\vee \xi+\int_t^Tg(s,\bar X_s,\bar Z_s){\rm ds}+\int_t^T{\rm d}V_s+n\int_t^T\left(\bar X_s-L_s\right)^-{\rm d}s\\
&&\Dis +\int_t^T\left(g^-(s,X_s,H_s){\rm d}s+{\rm d}C^-_s+{\rm d}V^-_s\right)-\int_t^T\bar Z_s{\rm d}B_s,\ \ t\in\T.
\end{array}
$$
Thus, by Corollary 3 again we know that $Y^n_t\leq \bar X_t$ for each $t\in \T$ and $n\geq 1$. That is to say, (C) holds true for $Y^n_\cdot$ and $g$.

Thus, we have verified that all conditions in Proposition 2 are satisfied with $g_n\equiv g$.
It then follows from Proposition 2 that (15) holds true. Furthermore, it follows from Corollary 3 that $Y_\cdot^n$ increases in $n$. Then we have
\begin{equation}
|Y^n_\cdot|\leq |Y^1_\cdot|+|\bar X_\cdot|\in \s^p.
\end{equation}
By combining (15) and (56) we can deduce that (18) holds true with
$$
\begin{array}{lll}
\eta &:=&\Dis C\left[\sup\limits_{t\in\T}
|Y^1_t|^p+|V|_T^p+\sup\limits_{t\in\T}|\bar X_t|^p+\left(\int_0^T f_t{\rm d}t\right)^p+1\right.\\
&&\ \ \ \ \ \Dis \left. +\left(\int_0^T |g(t,\bar X_t,0)|{\rm d}t\right)^p+\left(\int_0^T |g(t,0,0)|{\rm d}t\right)^p\right],
\end{array}
$$
where $C$ is a nonnegative constant depending on $p,\mu,\lambda,A,T$. \vspace{0.2cm}\hfill $\Box$

Making use of Theorems 1-2, Propositions 3, 5 and 7 together with Corollaries 2-3, we can prove the following existence and uniqueness results (see Theorems 4-5) for $L^p$ solutions of RBSDE (1) with $L^p\ (p>1)$ data under the assumptions of (H1), (H2) (resp. (H2w)), (H3), (H4w) (resp. (H4s)), (H5) and (H6), which also shows that the solutions can be approximated by the penalization method. \vspace{0.2cm}

{\bf Theorem 4}\ Let $p>1$, $V_\cdot\in\vcal^{p}$ and let $g$ satisfy (H1), (H2), (H3) and (H4w). Assume that (H5) and (H6) hold for $\xi$, $L_\cdot$ and some $X_\cdot$. For each $n\geq 1$, let $(Y_\cdot^n,Z_\cdot^n)\in \s^p\times\M^p$ be the unique solution of the penalization equation (16) with (17) (Recall Theorem 1). Then, there exists a triple $(Y_\cdot,Z_\cdot,K_\cdot)\in \s^p\times\M^p\times\vcal^{+,p}$ such that
\begin{equation}
\lim\limits_{n\To \infty}\left(\|Y_\cdot^n-Y_\cdot\|_{\s^p}+
\|Z_\cdot^n-Z_\cdot\|_{\M^p}+\left\| K^n_\cdot-K_\cdot\right\|_{\s^p}\right)=0.
\end{equation}
And, $(Y_\cdot,Z_\cdot,K_\cdot)$ is the unique solution of RBSDE $(\xi,g+{\rm d}V,L)$ in $\s^p\times\M^p\times\vcal^{+,p}$.\vspace{0.2cm}

{\bf Proof.}\  It follows from Proposition 5 or Corollary 1 that $Y_\cdot^n$ increases in $n$. By (i) and (ii) of Remark 2, we know that (H4) holds true, and (H2) and (H3) can imply (HH). At the same time, it follows from Proposition 7 that (18) holds. Thus, we have checked all conditions in Proposition 3. It then follows from Proposition 3 that there exists a solution $(Y_\cdot,Z_\cdot,K_\cdot)\in \s^p\times\M^p\times\vcal^{+,p}$ of RBSDE $(\xi,g+{\rm d}V,L)$ such that
\begin{equation}
\lim\limits_{n\To \infty}\left(\|Y_\cdot^n-Y_\cdot\|_{\s^p}+
\|Z_\cdot^n-Z_\cdot\|_{\M^p}\right)=0,
\end{equation}
and there exists a subsequence $\{K_\cdot^{n_j}\}$ of $\{K_\cdot^n\}$ such that
$\lim\limits_{j\To\infty}\sup\limits_{t\in\T}
|K_t^{n_j}-K_t|=0.$

In the sequel, we prove that
\begin{equation}
\lim\limits_{n\To\infty}\left\| \int_0^\cdot g(s,Y_s^n,Z_s^n){\rm d}s-\int_0^\cdot g(s,Y_s,Z_s){\rm d}s\right\|_{\s^p}=0.
\end{equation}
In fact, it follows from (H2) that $\as$, for each $n\geq 1$,
$$
\begin{array}{ll}
&\Dis |g(\cdot,Y_\cdot^n,Z_\cdot^n)-g(\cdot,Y_\cdot,Z_\cdot)|\\
\leq &\Dis  |g(\cdot,Y_\cdot^n,Z_\cdot^n)-g(\cdot,Y_\cdot^n,Z_\cdot)| +|g(\cdot,Y_\cdot^n,Z_\cdot)-g(\cdot,Y_\cdot,Z_\cdot)|\\
\leq &\Dis |g(\cdot,Y_\cdot^n,Z_\cdot)-g(\cdot,Y_\cdot,Z_\cdot)|
+\phi(|Z_\cdot^n-Z_\cdot|).
\end{array}
$$
Thus, making use of the following basic inequality (see \citet{Fan10} for details)
$$
\phi(x)\leq (m+2A)x+\phi\left({2A\over m+2A}\right),\ \ \RE\ x\geq 0,\ \ \RE m\geq 1
$$
together with H\"{o}lder's inequality, we get that for each $n,m\geq 1$,
\begin{equation}
\begin{array}{ll}
&\Dis \left\| \int_0^\cdot g(s,Y_s^n,Z_s^n){\rm d}s-\int_0^\cdot g(s,Y_s,Z_s){\rm d}s\right\|_{\s^p}\\
\leq & \Dis \left\| g(\cdot,Y_\cdot^n,Z_\cdot^n)
-g(\cdot,Y_\cdot,Z_\cdot)\right\|_{\hcal^p}\leq   \Dis \left\| g(\cdot,Y_\cdot^n,Z_\cdot)
-g(\cdot,Y_\cdot,Z_\cdot)\right\|_{\hcal^p}\\
&\Dis \ +(m+2A)T^{p-1\over p}\|Z_\cdot^n-Z_\cdot\|_{\M^p}+\phi({2A\over m+2A})T.
\end{array}
\end{equation}
Furthermore, note that $Y_t^1\leq Y_t^n\leq Y_t$ for each $n\geq 1$ and $t\in\T$. It follows from (H1) and (H2) together with Remark 1 that $\as$, for each $n\geq 1$, in view of (18), (iv) of Remark 2, and (4) in Lemma 3,
$$
\begin{array}{lll}
|g(\cdot,Y^n_\cdot,Z_\cdot)|&\leq & |g(\cdot,Y^1_\cdot,Z_\cdot)|
+|g(\cdot,Y_\cdot,Z_\cdot)|+2A|Y_\cdot
-Y^1_\cdot|+2A\\
&\leq & |g(\cdot,Y^1_\cdot,Z^1_\cdot)|
+|g(\cdot,Y_\cdot,Z_\cdot)|+2A|Y_\cdot
-Y^1_\cdot|+A|Z_\cdot
-Z^1_\cdot|+3A\in \hcal^p.
\end{array}
$$
Then, Lebesgue's dominated convergence theorem together with (H4w) and the fact of $\as, Y^n_\cdot\uparrow Y_\cdot$ yields that
\begin{equation}
\lim\limits_{n\To\infty}\left\| g(\cdot,Y_\cdot^n,Z_\cdot)
-g(\cdot,Y_\cdot,Z_\cdot)\right
\|_{\hcal^p}=0.
\end{equation}
Thus, letting first $n\To \infty$, and then $m\To\infty$ in (60), in view of (61), (58) and the fact that $\phi(\cdot)$ is continuous and $\phi(0)=0$, we get (59).

Finally, (57) follows from (58) and (59). And, the uniqueness part is a direct corollary of Proposition 5 or Corollary 1. The proof of Theorem 4 is completed.\vspace{0.2cm}\hfill $\Box$

{\bf Corollary 4}\ Let $p>1$, $V_\cdot^1, V_\cdot^2,\in\vcal^{p}$ and both $g^1$ and $g^2$ satisfy assumptions (H1), (H2), (H3) and (H4w). For $i=1,2$, assume that (H5) and (H6) hold for $\xi^i$, $L_\cdot^i$ and $X_\cdot^i$ associated with $g^i$, and that $(Y_\cdot^i,Z_\cdot^i,K_\cdot^i)\in \s^p\times \M^p\times \vcal^{+,p}$ is the unique solution of RBSDE $(\xi^i,g^i+{\rm d}V^i,L^i)$. If $\xi^1\leq \xi^2$, ${\rm d}V^1_\cdot\leq {\rm d}V^2_\cdot$, $L^1_\cdot=L^2_\cdot$, and
$$\as ,\ \ g^1(t,y,z)\leq g^2(t,y,z)$$
for each $(y,z)\in \R\times \R^d$, then $\ps$, ${d}K_t^1\geq {\rm d}K_t^2$ for each $t\in \T$.\vspace{0.2cm}

{\bf Proof.}\ For $n\geq 1$ and $i=1,2$, by Theorem 1 let $(Y_\cdot^{i,n},Z_\cdot^{i,n})\in \s^p\times\M^p$ be the unique solution of the following penalization BSDE:
$$
Y_t^{i,n}=\xi^i+\int_t^T g^i(s,Y_s^{i,n},Z_s^{i,n}){\rm d}s+\int_t^T {\rm d}V_s^i+\int_t^T {\rm d}K_s^{i,n}-\int_t^T Z_s^{i,n}{\rm d}B_s,\ \ t\in \T
$$
with
$$
K_t^{i,n}:=n\int_0^t\left(Y_s^{i,n}-L_s^i\right)^-{\rm d}s,\ \ t\in\T.\vspace{0.1cm}
$$
In view of the assumptions of Corollary 4, it follows from Proposition 5 that for each $n\geq 1$, $Y^{1,n}_\cdot\leq Y^{2,n}_\cdot$, and then
$$
K_{t_2}^{1,n}-K_{t_1}^{1,n}=n\int_{t_1}^{t_2}
\left(Y^{1,n}_s-L^1_s\right)^-{\rm d}s\geq n\int_{t_1}^{t_2}
\left(Y^{2,n}_s-L^2_s\right)^-{\rm d}s=K_{t_2}^{2,n}-K_{t_1}^{2,n}
$$
for every $n\geq 1$ and $0\leq t_1\leq t_2\leq T$. Since
$$
\|K^{1,n}_\cdot-K^1_\cdot\|_{\s^p}\To 0\ \ {\rm and}\ \ \|K^{2,n}_\cdot-K^2_\cdot\|_{\s^p}\To 0
$$
as $n\To\infty$ by Theorem 4, it follows that $\ps$,
$$K_{t_2}^1-K_{t_1}^1\geq K_{t_2}^2-K_{t_1}^2$$
for every $0\leq t_1\leq t_2\leq T$, which proves the desired result.\vspace{0.2cm}\hfill $\Box$

{\bf Remark 6}\ By (i) and (iii) of Remark 2, it is clear that Theorems 3-4 together with Corollary 4 strengthen the corresponding results for RBSDE (1) established in \citet{Kli12,Lep05} and \citet{Roz12}, which the stronger assumptions (H1s) and (H2s) than (H1) and (H2) are satisfied.\vspace{0.2cm}

{\bf Theorem 5}\ Let $p>1$, $V_\cdot\in\vcal^{p}$ and let $g$ satisfy (H1), (H2w), (H3) and (H4s). Assume that (H5) and (H6) hold for some $\xi$, $L_\cdot$ and $X_\cdot$. For each $n\geq 1$, let $(Y_\cdot^n,Z_\cdot^n)\in \s^p\times\M^p$ be the maximal (resp. minimal) solution of the penalization BSDE (16) with (17) (Recall Corollary 2). Then, there exists a solution $(Y_\cdot,Z_\cdot,K_\cdot)\in \s^p\times\M^p\times\vcal^{+,p}$ of RBSDE $(\xi,g+{\rm d}V,L)$ such that
$$
\lim\limits_{n\To \infty}\left(\|Y_\cdot^n-Y_\cdot\|_{\s^p}+
\|Z_\cdot^n-Z_\cdot\|_{\M^p}\right)=0,
$$
and there exists a subsequence $\{K_\cdot^{n_j}\}$ of $\{K_\cdot^n\}$ such that
$$\lim\limits_{j\To\infty}\sup\limits_{t\in\T}
|K_t^{n_j}-K_t|=0.$$

{\bf Proof.}\ It follows from Corollary 3 that $Y_\cdot^n$ increases in $n$. By (i) and (ii) of Remark 2 we know that (H2w) and (H3) imply (HH), and (H4) holds true. And, it follows from Proposition 7 that (18) holds true. Thus, we have checked all conditions in Proposition 3. Then the conclusion follows from Proposition 3. \vspace{0.2cm}\hfill $\Box$

{\bf Remark 7}\ Let us remark that it is not clear whether $(Y_\cdot,Z_\cdot)$ obtained in Theorem 5 is the maximal (resp. minimal) solution of RBSDE $(\xi,g+{\rm d}V,L)$ or not.\vspace{0.1cm}

The following Theorem 6 further proves that under the conditions of Theorem 5, RBSDE $(\xi,g+{\rm d}V,L_\cdot)$ admits both a minimal and a maximal solution in $\s^p\times\M^p\times \vcal^{+,p}$, which also shows that the solution can be approximated by some sequence of solutions of RBSDEs. The proof is based on Theorem 4, Corollaries 1 and 4, Propositions 1, 2, 4 and 5.\vspace{0.2cm}

{\bf Theorem 6}\ Let $p>1$, $V_\cdot\in\vcal^{p}$ and let $g$ satisfy (H1), (H2w), (H3) and (H4s). Assume that (H5) and (H6) hold for some $\xi$, $L_\cdot$ and $X_\cdot$. Then, RBSDE $(\xi,g+{\rm d}V,L_\cdot)$ admits a minimal (resp. maximal) solution $(Y_\cdot,Z_\cdot, K_\cdot)$ in $\s^p\times\M^p\times \vcal^{+,p}$. i.e., if $(Y'_\cdot,Z'_\cdot,K'_\cdot)$ is also a solution of RBSDE $(\xi,g+{\rm d}V, L_\cdot)$ in the space $\s^p\times\M^p\times \vcal^{+,p}$, then $\ps$, $Y_t\leq Y'_t$ (resp. $Y_t\geq Y'_t$) for each $t\in \T$.\vspace{0.2cm}

{\bf Proof.}\ Assume that $p>1$, $V_\cdot\in\vcal^{p}$ and the generator $g$ satisfies (H1) with $\rho(\cdot)$, (H2w) with $f_\cdot$, $\mu$ and $\lambda$, (H3) with $\varphi_\cdot(r)$ and (H4s). Assume further that (H5) and (H6) hold for some $\xi$, $L_\cdot$ and $X_\cdot$.

We first show the existence of the minimal solution. In view of the assumptions of $g$, it is not very hard to prove that for each $n\geq 1$ and $(y,z)\in \R\times\R^d$, the following function
\begin{equation}
g_n(\omega,t,y,z):=\inf\limits_{u\in\R^d}
\left[g(\omega,t,y,u)+(n+2\lambda)|u-z|\right]
\end{equation}
is well defined and $(\F_t)$-progressively measurable, $\as$, $g_n$ increases in $n$ and  converges locally uniformly in $(y,z)$ to the generator $g$ as $n\To \infty$, and all $g_n$ satisfy (H1) with the same $\rho(\cdot)$, (H2s) with $n+2\lambda$, (H3) with the same $\varphi_\cdot(r)+\mu r+2f_\cdot$ and (H4). In addition, we can also prove that for each $n\geq 1$, $g_n$ and $g$ satisfy (14) (appearing in Proposition 2) with the same $f_\cdot$, $\mu$ and $\lambda$, and then $\as$, for each $n\geq 1$ and each $(y,z)\in \R\times\R^d$,
$$
|g_n(\cdot,y,z)|\leq |g(\cdot,y,0)|+f_\cdot+\mu |y|+\lambda |z|\leq |g(\cdot,0,0)|+f_\cdot+\varphi_\cdot(|y|)+\mu |y|+\lambda |z|.
$$
That is to say, all $g_n$ satisfy (HH) with the same parameters.

Note that (H2s) implies (H2), and (H4) implies (H4w) by (i) of Remark 2. It then follows from Theorem 4 that there exists a unique solution $(Y_\cdot^n,Z_\cdot^n, K_\cdot^n)\in \s^p\times\M^p\times \vcal^{+,p}$ of RBSDE $(\xi,g_n+{\rm d}V,L)$ for each $n\geq 1$. By Corollary 1, $Y_\cdot^n$ increases in $n$. Furthermore, it follows from (14) and (H6) that for each $n\geq 1$,
$$
|g_n(\cdot,X_\cdot,0)|\leq |g(\cdot,X_\cdot,0)|+f_\cdot+\mu |X_\cdot|\in \hcal^p.
$$
Then, by Corollary 4, $K_\cdot^n$ decreases in $n$.

In the sequel, we show that (C) appearing in Proposition 2 holds true for $Y^n_\cdot$ and $g$. In fact,
let
$$
\underline{g}(\cdot,y,z):=g(\cdot,y,0)-f_\cdot-\mu |y|-\lambda |z|
$$
and
$$
\bar g(\cdot,y,z):=g(\cdot,y,0)+f_\cdot+\mu |y|+\lambda |z|.\vspace{0.2cm}
$$
Then both $\underline{g}$ and $\bar g$ satisfy (H1), (H2s), (H3) and (H4s),
$$
\underline{g}(\cdot, X_\cdot,0)=g(\cdot, X_\cdot,0)-f_\cdot-\mu |X_\cdot|\in \hcal^p,
$$
$$
\bar g(\cdot, X_\cdot,0)=g(\cdot, X_\cdot,0)+f_\cdot+\mu |X_\cdot|\in \hcal^p,\vspace{0.1cm}
$$
and by (14) for each $n\geq 1$,
$$\underline{g}\leq g_n\leq \bar g.\vspace{0.1cm}$$
Thus, it follows from Theorem 4 that RBSDE $(\xi,\underline{g}+{\rm d}V,L)$ and RBSDE $(\xi,\bar g+{\rm d}V,L)$ admit respectively a unique solution $(\underline{Y}_\cdot,\underline{Z}_\cdot, \underline{K}_\cdot)\in \s^p\times\M^p\times \vcal^{+,p}$ and $(\bar Y_\cdot,\bar Z_\cdot, \bar K_\cdot)\in \s^p\times\M^p\times \vcal^{+,p}$, and by Corollary 1, we know that $\ps$,
\begin{equation}
\underline{Y}_t\leq Y^n_t\leq \bar Y_t
\end{equation}
for each $t\in \T$ and $n\geq 1$. In addition, in view of (i) and (ii) of Remark 2, by Proposition 1 with $g_n\equiv\bar g$ we know that $\bar g(\cdot, \bar Y_\cdot,\bar Z_\cdot)\in\hcal^p$, and then
\begin{equation}
g(\cdot, \bar Y_\cdot,0)=\bar g(\cdot, \bar Y_\cdot,\bar Z_\cdot)-f_\cdot-\mu|\bar Y_\cdot|-\lambda|\bar Z_\cdot|\in\hcal^p.\vspace{0.1cm}
\end{equation}
By (63) and (64) we know that (C) is true for $Y^n_\cdot$ and $g$.

Now we have checked that all conditions in Proposition 2 are satisfied. It then follows from Proposition 2 that (15) holds true. Furthermore, in view of (63), we can deduce that (37) appearing in Proposition 4 holds true with
$$
\begin{array}{lll}
\eta &:=&\Dis C\left[\sup\limits_{t\in\T}
|\underline{Y}_t|^p+|V|_T^p+\sup\limits_{t\in\T}|\bar Y_t|^p+\left(\int_0^T f_t{\rm d}t\right)^p+1\right.\\
&&\ \ \ \ \ \Dis \left. +\left(\int_0^T |g(t,\bar Y_t,0)|{\rm d}t\right)^p+\left(\int_0^T |g(t,0,0)|{\rm d}t\right)^p\right],
\end{array}
$$
where $C$ is a nonnegative constant depending on $p,\mu,\lambda,A,T$. Hence, all conditions in Proposition 4 are satisfied, and then it follows from Proposition 4 that RBSDE $(\xi,g+{\rm d}V,L)$ admits a solution $(Y_\cdot,Z_\cdot,K_\cdot)\in \s^p\times\M^p\times\hcal^p$ such that
\begin{equation}
\lim\limits_{n\To\infty}(\|Y_\cdot^n
-Y_\cdot\|_{\s^p}+\|Z_\cdot^n
-Z_\cdot\|_{\M^p}+\|K_\cdot^n
-K_\cdot\|_{\s^p})=0.\vspace{-0.2cm}
\end{equation}

Finally, let us show that $(Y_\cdot,Z_\cdot,K_\cdot)$ is just the minimal solution of RBSDE $(\xi,g+{\rm d}V,L)$ in $\s^p\times\M^p\times\hcal^p$. In fact, if $(Y'_\cdot,Z'_\cdot,K'_\cdot)$ is also a solution of RBSDE $(\xi,g+{\rm d}V,L)$ in $\s^p\times\M^p\times\hcal^p$, then noticing that for each $n\geq 1$, $g_n\leq g$ and $g_n$ satisfies (H1) and (H2) due to (H2s)$\Rightarrow$(H2), it follows from Proposition 5 that $\ps$, $Y_t^n\leq Y'_t$ for each $t\in\T$ and $n\geq 1$. Thus, by (65), $\ps$, $Y_t\leq Y'_t$ for each $t\in \T$.

As for the case of the maximal solution, we only need to replace (62) with (53). And, by a similar argument as above we can obtain the desired result. The proof of Theorem 6 is then completed. \vspace{0.2cm}\hfill$\Box$

{\bf Remark 8}\ It follows from (i) and (iii) of Remark 2 that Theorem 6 improves Theorem 5.1 in \citet{Xu08}, where the stronger assumption (H1s) than (H1) is satisfied, the barrier $L_\cdot$ is assumed to be bounded and only $L^2$ solution is considered.\vspace{0.2cm}

By Corollaries 1, 4 and the proof of Theorem 6 it is not hard to verify that under assumptions (H1), (H2w), (H3), (H4s), (H5) and (H6), the comparison theorem for the maximal (resp. minimal) $L^p$ solutions of RBSDEs with $L^p\ (p>1)$ data is true. More precisely, we have\vspace{0.1cm}

{\bf Proposition 8}\ Let $p>1$ and for $i=1,2$, assume that $V_\cdot^i\in \vcal^p$, $g^i$ satisfies (H1), (H2w), (H3) and (H4s), $\xi^i$, $L_\cdot^i$ and $X_\cdot^i$ satisfy (H5) and (H6) associated with $g^i$, and $(Y_\cdot^i,Z_\cdot^i,K_\cdot^i)\in \s^p\times \M^p\times \vcal^{+,p}$ is the maximal (resp. minimal) solution of RBSDE $(\xi^i,g^i+{\rm d}V^i,L^i)$ (Recall Theorem 6). If $\xi^1\leq \xi^2$, ${\rm d}V^1_\cdot\leq {\rm d}V^2_\cdot$, $L^1_\cdot\leq L^2_\cdot$, and
$$\as,\ \ g^1(t,y,z)\leq g^2(t,y,z)$$
for each $(y,z)\in \R\times \R^d$, then $\ps$, $Y_t^1\leq Y_t^2$ for each $t\in \T$. Furthermore, if $L^1_\cdot=L^2_\cdot$, then $\ps$, ${\rm d}K_t^1\geq {\rm d}K_t^2$ for each $t\in \T$.\vspace{0.2cm}

{\bf Remark 9}\ In a subsequent work we will further discuss the problem on the existence and uniqueness for $L^p$ solutions of RBSDE (1) with $L^p\ (p>1)$ or $L^1$ data, irregular barriers and generators satisfying (H1), (H2) or (H2w), and (H3). \vspace{0.6cm}




\begin{thebibliography}{00}

\setlength{\baselineskip}{13.5pt}
\setlength{\itemsep}{1.8mm}

\bibitem[Aman(2009)]{Ama09}Aman, A. (2009). $L^p$-solutions of reflected generalized BSDEs with non-Lipschitz coefficients. {\it Random Oper. Stoch. Equ.} {\bf 17}, 201-219.

\bibitem[Bayraktar and Yao(2012)]{Bay12} Bayraktar, E., Yao, S. (2012). Quadratic reflected BSDEs with unbounded obstacles. {\it Stochastic Process. Appl.} {\bf 122}, 1155-1203.

\bibitem[Bayraktar and Yao(2014)]{Bay14} Bayraktar, E., Yao, S. (2014). Doubly reflected BSDEs with integrable parameters and related Dynkin games. arXiv: 1412.2053v1 [math.PR].

\bibitem[Briand, Delyon, Hu, Pardoux and Stoica(2003)]{Bri03}Briand, Ph., Delyon, B., Hu, Y., Pardoux, E., Stoica, L. (2003). $L^p$ solutions of backward stochastic differential equations. {\it Stochastic Process. Appl.} {\bf 108}, 109-129.

\bibitem[Briand, Lepeltier and San Martin(2007)]{Bri07}Briand, Ph., Lepeltier, J.-P., San Martin, J. (2007). One-dimensional BSDEs whose coefficient is monotonic in $y$ and non-Lipschitz in $z$. {\it Bernoulli} {\bf 13}, 80-91.


\bibitem[El Karoui, Kapoudjian, Pardoux, Peng and Quenez(1997)]{El97}El Karoui, N., Kapoudjian, C., Pardoux, E., Peng, S., Quenez, M.-C. (1997). Reflected solutions of backward SDEs, and related obstacle problems for PDE's. {\it Ann. Probab.} {\bf 25}, 702-737.

\bibitem[El Karoui, Pardoux and Quenez(1997)]{El97a}El Karoui, N., Pardoux, E., Quenez, M.-C. (1997). Reflected backward SDEs and American options, in {\it Numerical Methods in Finance}, eds. Robers L. and Talay D. (Cambridge Univ. Press), pp. 215-231.

\bibitem[El Karoui, Peng and Quenez(1997)]{El97b}El Karoui, N., Peng, S., Quenez, M.C. (1997). Backward stochastic differential equations in finance. {\it Math. Finance} {\bf 7}, 1-72.

\bibitem[Fan(2015)]{Fan15}Fan, S. (2015). $L^p$ solutions of multidimensional BSDEs with weak monotonicity and general growth generators. {\it Journal of Mathematical Analysis and Applications} {\bf 432}, 156-178.

\bibitem[Fan and Jiang(2010)]{Fan10}Fan, S., Jiang, L. (2010). Uniqueness result for the BSDE whose generator is monotonic in $y$ and uniformly continuous in $z$. {\it C. R. Acad. Sci. Paris, Ser. I} {\bf 348}(1-2), 89-92.

\bibitem[Fan and Jiang(2012)]{Fan12}Fan, S., Jiang, L. (2012). A generalized comparison theorem for BSDEs and its Applications. {\it J. Theor. Probab.} {\bf 25}, 50-61.

\bibitem[Fan, Jiang and Davison(2013)]{Fan13}Fan, S., Jiang, L., Davison, M.(2013). Existence and uniqueness result for multidimensional BSDEs with generators of Osgood type. {\it Front. Math. China}, {\bf 8}(4), 811-824.

\bibitem[Hamad\`{e}ne(2002)]{Ham02}Hamad\`{e}ne, S. (2002). Reflected BSDE's with discontinuous barrier and application. {\it Stochastics and Stochastics Reports} {\bf 74}(3-4), 571-596.

\bibitem[Hamad\`{e}ne and Lepeltier(2000)]{Ham00}Hamad\`{e}ne, S., Lepeltier, J.-P. (2000). Reflected BSDEs and mixed game problem. {\it Stochastic Process. Appl.} {\bf 85}, 177-188.

\bibitem[Hamad\`{e}ne and Popier(2012)]{Ham12}Hamad\`{e}ne, S., Popier, A. (2012). $L^p$-solutions for reflected backward stochastic differential equations. {\it Stoch. Dyn.} {\bf 12}, 1150016, 35 pp.

\bibitem[Hamad\`{e}ne and Zhang(2010)]{Ham10}Hamad\`{e}ne, S., Zhang, J. (2010). Switching problem and related system of reflected backward SDEs. {\it Stochastic Process. Appl.} {\bf 120}(4), 403-426.

\bibitem[Hu and Tang(2010)]{Hu08}Hu, Y., Tang, S. (2010). Multi-dimensional BSDE with oblique Reflection and optimal switching. {\it Probab. Theory and Related Fields} {\bf 147}(1-2), 89-121.

\bibitem[Hu and Tang(2014)]{Hu14}Hu, Y., Tang, S.(2014). Multidimensional backward stochastic differential equations of diagonally quadratic generators. {\it arXiv: 1408.4579v1 [math.PR]}.

\bibitem[Hua, Jiang and Shi(2013)]{Hua13}Hua, W., Jiang, L., Shi, X. (2013). Infinite time interval RBSDEs with non-Lipschitz coefficients. {\it Journal of the Korean Statistical Society} {\bf 42}, 247-256.

\bibitem[Izumi(2013)]{Izu13}Izumi, Y. (2013). The $L^p$ Cauchy sequence for one-dimensional BSDEs with linear growth generators {\it Statist. Probab. Lett.} {\bf 83}, 1588-1594.

\bibitem[Jia(2008)]{Jia08}Jia, G. (2008). A uniqueness theorem for the solution of backward stochastic differential equations. {\it C. R. Acad. Sci. Paris, Ser. I} {\bf 346}, 439-444.

\bibitem[Jia(2010)]{Jia10}Jia, G. (2010). Backward Stochastic differential equations with a uniformly continuous generator and related $g$-expectation. {\it Stochastic Process. Appl.} {\bf 120}(11), 2241-2257.

\bibitem[Jia and Xu(2008)]{Jia08a}Jia, G., Xu, M. (2008). Construction and Uniqueness for reflected BSDE under linear increasing condition. {\it arXiv:0801.3718v1 [math.SG]}.

\bibitem[Klimsiak(2012)]{Kli12}Klimsiak, T. (2012). Reflected BSDEs with monotone generator. {\it Electron. J. Probab.} {\bf 107}, 1-25.

\bibitem[Klimsiak(2013)]{Kli13}Klimsiak, T. (2013). BSDEs with monotone generator and two irregular reflecting barriers. {\it Bull. Sci. math.} {\bf 137}, 268-321.

\bibitem[Kobylanski, Lepeltier, Quenez and Torres (2002)]{Kob02}Kobylanski, M., Lepeltier, J.-P., Quenez, M.C., Torres, S. (2002). Reflected BSDE with superlinear quadratic coefficient. {\it Probab. Math. Statist.} {\bf 22}, 51-83.

\bibitem[Lepeltier, Matoussi and Xu (2005)]{Lep05}Lepeltier, J.-P., Matoussi, A., Xu, M. (2005). Reflected bakcward stochastic differential equations under monotonicity and general increasing growth conditions. {\it Adv. in Appl. Probab.} {\bf 37}, 134-159.

\bibitem[Lepeltier and Xu (2005)]{Lep05a}Lepeltier, J.-P., Xu, M. (2005). Penalization method for reflected backward stochastic differenctial equations with one R.C.L.L. barrier.  {\it Statist. Probab. Lett.} {\bf 75}, 58-66.

\bibitem[Ma, Fan and Song(2013)]{Ma13}Ma, M., Fan, S., Song, X.(2013). $L^p (p >1)$ solutions of backward stochastic differential equations with monotonic and uniformly continuous generators. {\it Bull. Sci. math} {\bf 137}, 97-106.

\bibitem[Mao(1995)]{Mao95}Mao, X. (1995). Adapted solutions of backward stochastic differential equations with non-Lipschitz cofficients. {\it Stochastic Process. Appl.} {\bf 58}, 281-292.

\bibitem[Matoussi(1997)]{Mat97}Matoussi, A. (1997). Reflected solutions of backward stochastic differential equations with continuous coefficient. {\it Statist. Probab. Lett.} {\bf 34}, 347-354.

\bibitem[Pardoux(1999)]{Par99}Pardoux, E. (1999). BSDEs, weak convergence and homogenization of semilinear PDEs. Nonlinear Analysis, Differential Equations and Control (Montreal, QC,1998). Kluwer Academic Publishers, Dordrecht, pp.503-549.

\bibitem[Pardoux and Peng(1990)]{Par90}Pardoux, E., Peng, S. (1990). Adapted solution of a backward stochastic differential equation. {\it Systems Control Lett.} {\bf 14}, 55-61.

\bibitem[Peng(1999)]{Peng99}Peng, S. (1999). Monotonic limit theorem of BSDE and nonlinear decomposition theorem of Doob-Meyer's type. {\it Probab. Theory Related Fields} {\bf 113},473-499.

\bibitem[Peng and Xu(2005)]{Peng05}Peng, S., Xu, M. (2005). The smallest $g$-supermartingale and reflected BSDE with single and double $L^2$ obstacles. {\it Ann. Inst. H. Poincar\'{e} Probab. Statist.} {\bf 41}, 605-630.

\bibitem[Peng and Xu(2010)]{Peng10}Peng, S., Xu, M. (2010). Reflected BSDE with a constraint and its applications in an incomplete market. {\it Bernoulli} {\bf 16}(3), 614-640.

\bibitem[Rozkosz and S{\l}omi\'{n}ski(2012)]
    {Roz12}Rozkosz, A., S{\l}omi\'{n}ski, L. (2012). $L^p$ solutions of reflected BSDEs under monotonicity condition. {\it Stochastic Process. Appl.} {\bf 122}, 3875-3900.

\bibitem[Xu(2008)]{Xu08}Xu, M. (2008). Backward stochastic differential equations with reflection and weak assumptions on the coefficients. {\it Stochastic Process. Appl.} {\bf 118}, 968-980.


\end{thebibliography}
\end{document}